# Signature cumulants, ordered partitions, and independence of stochastic processes

PATRIC BONNIER[*] and HARALD OBERHAUSER[†]

*Mathematical Institute, Andrew Wiles Building, Radcliffe Observatory Quarter, Woodstock Road, OX2 6GG, Oxford, United Kingdom. E-mail:* [*]*Patric.Bonnier@maths.ox.ac.uk;* [†]*Harald.Oberhauser@maths.ox.ac.uk*

The sequence of so-called signature moments describes the laws of many stochastic processes in analogy with how the sequence of moments describes the laws of vector-valued random variables. However, even for vector-valued random variables, the sequence of cumulants is much better suited for many tasks than the sequence of moments. This motivates us to study so-called signature cumulants. To do so, we develop an elementary combinatorial approach and show that in the same way that cumulants relate to the lattice of partitions, signature cumulants relate to the lattice of so-called "ordered partitions". We use this to give a new characterisation of independence of multivariate stochastic processes. Finally, we construct a family of unbiased minimum-variance estimators of signature cumulants and show that even for the simple example of a diffusion with constant drift and volatility, such signature cumulant estimators outperform signature moment estimators.

*Keywords:* cumulants; geometric rough paths; partitions; path signatures; stochastic processes

## 1. Introduction

The sequence of moments $(\mu_X^m)_{m\geq 1}$ of an $\mathbb{R}^d$-valued random variable $X$,

$$\left(\mu_X^m\right)_{m\geq 1} \quad \text{where } \mu_X^m := \mathbb{E}\left[X^{\otimes m}\right] \in \left(\mathbb{R}^d\right)^{\otimes m}, \tag{1.1}$$

plays a fundamental role in statistics and probability theory since it captures statistical properties of the law of $X$ in a graded manner. An undesirable property of moments – especially when only a finite number of samples from $X$ is available to estimate them – is that lower order moments can dominate higher order moments. For example, even for $m = 2$ and $d = 1$

$$\mu_X^2 = \left(\mu_X^1\right)^2 + \text{Var}(X),$$

and the variance $\text{Var}(X) := \mathbb{E}[(X - \mu_X^1)^2]$ is a much more sensible "order $m = 2$" statistic than $\mu_X^2$. This method of subtracting lower order moments from higher order moments to obtain better statistics can be continued beyond level $m = 2$; for $m = 3$ one gets the classical notion of "skewness", for $m = 4$ the "tailedness", etc. In general, the resulting sequence is known as the *cumulants*, which differ from the centralised moments in that one allows for subtractions of any lower order moment. Another advantage cumulants have over moments is that many probabilistic properties are often easier expressed in cumulants than in moments.

**Example 1.1.** One good reason to use cumulants instead of moments is that their sample statistics typically have lower variance. To illustrate this on a elementary example assume that $X \sim N(\mu, \sigma^2)$ follows a normal distribution. Given $N$ independent samples $X_1, \ldots, X_N$ of $X$ the minimum variance unbiased estimators for its second moment and cumulant are respectively,

$$\hat{\mu}_X^2 = \frac{1}{N} \sum_{i=1}^N X_i^2, \qquad \hat{\kappa}_X^2 = \frac{1}{N-1} \sum_{i=1}^N \left( X_i - \frac{1}{N} \sum_{j=1}^N X_j \right)^2.$$

By explicit computation in the same fashion as outlined in Appendix C one sees that

$$\mathrm{Var}(\hat{\mu}_X^2) = \mathrm{Var}(\hat{\kappa}_X^2) + \frac{2}{N}\left( \mu^4 - \mu^2 \sigma^2 - \frac{2\sigma^4}{N-1} \right).$$

It follows that when $\mu$ is large compared to $\sigma^2$, it is much preferable to use the cumulant estimator rather than the moment estimator.

Many real-world observations have an inherent sequential structure and one ultimately deals with path-valued samples (i.e., sample paths of a stochastic process) rather than vector-valued samples. For a path-valued random variable, the sequence of so-called signature moments

$$(\boldsymbol{\mu}_X^m)_{m \geq 0} = \left( \mathbb{E}\left[ \int dX^{\otimes m} \right] \right)_{m \geq 0}$$

is a natural replacement for the classical moment sequence (1.1). Indeed, motivated by insights from stochastic analysis [5], these signature moments have recently found applications in statistics and machine learning such as parameter estimation for stochastic differential equations [16] or hypothesis testing for laws of stochastic processes [6]. Motivated by the nice statistical properties of cumulants over moments in the case of vector-valued observations, the aim of this article is to study such properties for signature cumulants and to derive estimators for such signature cumulants. For example, in view of Example 1.1 one can ask whether cumulants can yield better estimators for estimating diffusion and drift in an SDE and we invite the reader to have a glance at Example 4.3 for an application of our results to a simple SDE with constant drift and volatility. Similarly, signature cumulants between stochastic processes vanish if and only if the stochastic processes are independent which allows for independence testing of time series.

### 1.1. Poset of partitions, moments and cumulants

Before further discussing the path-valued case of signature cumulants, we briefly introduce notation and recall results for classical cumulants below. The following result is classical, see [14] for history and overview. We spell it out in notation unusual in statistics however, as this formulation generalises naturally to our extension. For the rest of the paper, we will denote by $V$ a $d$-dimensional vector space with basis $\{e_1, \ldots, e_d\}$. The *tensor algebra* and its topological dual the *extended tensor algebra* of $V$ are defined as

$$\mathrm{T}(V) := \bigoplus_{m \geq 0} V^{\otimes m}, \qquad \mathrm{T}((V)) := \prod_{m \geq 0} V^{\otimes m}.$$

The pairing $\langle \cdot, \cdot \rangle$ between T($V$) and T(($V$)) is defined as $\langle s, t \rangle = \sum s_\tau t_\tau$ for $s = \sum s_\tau e_\tau$, $t = \sum t_\tau e_\tau$ where all sums are taken over multi-indices $\tau \in \bigcup_m \{1, \ldots, d\}^m$ and $e_\tau := e_{i_1} \otimes \cdots \otimes e_{i_m}$.

**Theorem 1.1 (Classical cumulants).** *Let $X$ be a $V$-valued random variable such that that the sequence of moments*

$$\mu_X := (\mu_X^m)_{m \geq 0} := (\mathbb{E}[X^{\otimes m}])_{m \geq 0} \in \mathrm{T}((V))$$

*is well-defined and characterises the law of $X$. We call*[1]

$$\kappa_X := (\kappa_X^m)_{m \geq 0} := \pi_{\mathrm{Sym}}(\log \mathbb{E}[\exp(X)]) \in \mathrm{T}((V)) \quad (1.2)$$

*the cumulant of $X$. For $\tau = (i_1, \ldots, i_m) \in \{1, \ldots, d\}^m$ it holds that*

1. *(compensated moments)*

$$\langle \kappa_X, e_\tau \rangle = \sum_{\mathbf{a} \in \mathcal{P}(\tau)} (-1)^{|\mathbf{a}|-1} (|\mathbf{a}| - 1)! \prod_{i=1}^k \langle \mu_X, e_{a_i} \rangle,$$

   *where the sum ranges over all partitions $\mathbf{a} = \{a_1, \ldots, a_k\}$ of $\tau$.*

2. *(moment bijection) There exist a bijection $\kappa_X \mapsto \mu_X$ explicitly given as*

$$\langle \mu_X, e_\tau \rangle = \sum_{\mathbf{a} \in \mathcal{P}(\tau)} \prod_{i=1}^k \langle \kappa_X, e_{a_i} \rangle,$$

   *where again the sum ranges over all partitions of $\tau$.*

3. *(independence) For $I, J \subset \{1, \ldots, d\}$ the families $(\langle X, e_i \rangle)_{i \in I}$ and $(\langle X, e_j \rangle)_{j \in J}$ are independent if and only if*

$$\langle \kappa_X, e_{\tau_1} \otimes e_{\tau_2} \rangle = 0 \quad \textit{for any } \tau_1 \in I^\star, \tau_2 \in J^\star,$$

*where $I^\star := \bigcup_m I^m$ and $J^\star$ is defined similarly.*

Item 1 connects cumulants with subtraction of lower moments from higher moments – thus also the term "compensated moments" – with definition (1.2) motivated by the Laplace transform. Item 2 implies that, like the sequence of moments, the sequence of cumulants characterises the law of $X$. Finally, the third item demonstrates that independence of random variables is equivalent to their cross-cumulants vanishing.

In two seminal papers [19,20], T. Speed showed that results such as the above can be proven in short and elegant fashion using the lattice of partitions on $\{1, \ldots, d\}$. The aim of this paper is to replace probability measures on $\mathbb{R}^d$ by probability measures on spaces of paths and to develop such a combinatorial approach and apply it to derive efficient estimators. As we will see, so-called signature cumulants of laws of stochastic processes are intimately linked to the lattice of "ordered partitions" which has a rich combinatorial structure.

---

[1] Using the notation $1 + x_1 + x_2 + \cdots$ for $(x_m)_{m \geq 0}$, the exponential and logarithm are defined as $\exp(x) = \sum_{m \geq 0} \frac{1}{m!} x^{\otimes m}$, $\log(1 + x) = \sum_{m \geq 1} \frac{(-1)^{m-1}}{m} x^{\otimes m}$. We denote with $\pi_{\mathrm{Sym}}(v_1 \otimes \cdots \otimes v_m) = \sum_\sigma v_{\sigma(1)} \otimes \cdots \otimes v_{\sigma(m)}$ the projection onto the symmetric tensors, save the $\frac{1}{m!}$ normalising factor.

## The poset of ordered partitions, signature moments and cumulants

Our results apply to any stochastic process for which a "rough path lift" exists, such as continuous semi-martingales. However, even for stochastic processes with smooth sample paths, it is interesting to spell out our "path-space version" of Theorem 1.1 which appears below as motivation for the rest of the paper. To speak of moments of measures on paths, we first need a notion of monomials of paths. If the path is smooth, this role is taken by iterated integrals. This sequence of tensors is called the *path signature*.

**Definition 1.1.** Define for $x \in \mathrm{BV}(V) := \bigcup_{T>0} \{x \in C([0, T], V) : \|x\|_{1\text{-var}} < \infty\}$ its signature $S(x) = (S_m(x))_{m \geq 0} = (\int dx^{\otimes m})_{m \geq 0}$ as

$$\int dx^{\otimes 0} := 1 \quad \text{and} \quad \int dx^{\otimes m} := \int_{0 < t_1 < \cdots < t_m < T} dx(t_1) \otimes \cdots \otimes dx(t_m) \quad \text{for } m \geq 1.$$

The study of how such iterated integrals reflect properties of the path $x$ goes back to Chen [4], see also [10]. If one applies this to random paths $X$, this naturally leads to the notion of signature moments $\mathbb{E}[S(X)] \in T(\!(V)\!)$. More recently such moments have found applications in statistics, starting with [16] where a "method of signature moments" is derived in analogy to the classical method of moments. We state Theorem 1.2 to give the reader an idea of the kind of statements we are interested in by using notation defined in the subsequent Section 2 and Section 3; but even without the concrete definitions of the poset $\mathrm{Orp}(w_1, \ldots, w_k)$, the shuffle ш, the unparametrised paths $\mathcal{R}_1(V)$, and the generalised signature moments $\mu_X(\mathbf{a})$ and cumulants $\kappa_X(\mathbf{a})$ at hand, the connection between Theorem 1.2 and Theorem 1.1 is evident. In fact, as we will in see in Example 3.3, Theorem 1.1 follows from Theorem 1.2 by tensor symmetrisation.

**Theorem 1.2 (Signature cumulants).** *Let $X$ be a $\mathcal{R}_1(V)$-valued random variable, such that the signature moments*

$$\boldsymbol{\mu}_X := \left(\mu_X^m\right)_{m \geq 0} := \left(\mathbb{E}[S_m(X)]\right)_{m \geq 0} \in T(\!(V)\!)$$

*are well-defined and characterise the law of $X$. We call*

$$\boldsymbol{\kappa}_X := \left(\kappa_X^m\right)_{m \geq 0} := \log \mathbb{E}[S(X)] \in T(\!(V)\!)$$

*the* signature cumulant *of $X$. For $\tau, \tau_1, \ldots, \tau_k \in \bigcup_m \{1, \ldots, d\}^m$ it holds that*

1. (compensated moments)

$$\langle \boldsymbol{\kappa}_X, e_{\tau_1} \text{ш} \cdots \text{ш} e_{\tau_k} \rangle = \sum_{\mathbf{a} \in \mathrm{Orp}(\tau_1, \ldots, \tau_k)} (-1)^{|\mathbf{a}|-1} \frac{\mathbf{a}!}{|\mathbf{a}|} \mu_X(\mathbf{a}),$$

*where the sum is taken over all ordered partitions of the tuple $(\tau_1, \ldots, \tau_k)$,*

2. (moment bijection)

$$\langle \boldsymbol{\mu}_X, e_\tau \rangle = \sum_{\mathbf{a} \in \mathrm{Orp}(\tau)} \frac{1}{|\mathbf{a}|!} \kappa_X(\mathbf{a}),$$

where the sum is taken over all ordered partitions of $\tau$,

3. (independence) *For* $I, J \subset \{1, \ldots, d\}$ *the families* $(\langle X, e_i \rangle)_{i \in I}$ *and* $(\langle X, e_j \rangle)_{j \in J}$ *are independent[2] if and only if*

$$\langle \kappa_X, e_{i^\star} \sqcup e_{j^\star} \rangle = 0 \quad \textit{for any } i^\star \in I^\star, j^\star \in J^\star.$$

We emphasise that as for classical cumulants, such results can in principle be derived by direct but lengthy calculations that compare coefficients; or, in contrast, an algebraic approach that relies on the well developed theory of Lie polynomials and Hall bases [18]; or even a Hopf algebraic perspective that has recently been applied to free cumulants [7] as well as Wick polynomials [8]. All of these approaches have their own merits, but as we hope to demonstrate, the combinatorial approach we present here via "ordered partitions" is intuitive, leads to simple proofs, and allows for a self-contained treatment throughout. Möbius inversion and basic properties of U-statistics are the only results we use for which we don't give full proofs. In short, the lattice of ordered partitions is useful for signature cumulants for the same reasons that the lattice of partitions is useful for cumulants.

**Remark 1.1.** There are other notions of non-commutative cumulants such as the Boolean, monotone, and free cumulants that arise in non-commutative probability theory. For instance, the free cumulants are studied using "non-crossing partitions", see [11,15,21]. While similar in spirit, our "ordered partitions" have a different combinatorial structure. Yet another related research area is the use of cumulants in Wiener Chaos expansions, see [17] for a survey, but the combinatorics are given by the standard poset of partitions.

## Structure

- In Section 2, we give some combinatorial background and define the set of ordered partitions. We then go on to show some of their properties in the setting that applies to us.
- In Section 3, we give some background on rough paths and study their signature cumulants.
- In Section 4, we construct an unbiased minimum variance estimator for signature cumulants and show that already for the simple case of a diffusion with constant drift and volatility, signature cumulants provide more efficient estimators than signature moments.

Appendix A contains technical details of independence of rough paths, Appendix B contains technical information about Tree-like equivalence of paths, and Appendix C contains calculations for Example 4.3.

---

[2]As usual we call two stochastic processes $X$ and $Y$ independent if their $\sigma$-algebras are independent. We emphasize that this is much stronger than instantaneous independence, i.e that $X_t$ and $Y_t$ are independent for every $t \geq 0$.

| Symbol | Meaning | Page |
|---|---|---|
| | Tensors | |
| $V$ | a $d$-dimensional vector space | 2728 |
| T($V$) | $\bigoplus V^{\otimes m}$ the free algebra on the vector space $V$ | 2728 |
| T(($V$)) | $\prod V^{\otimes m}$ the dual of the free algebra T($V$) | 2728 |
| $\langle \cdot, \cdot \rangle$ | the pairing between T($V$) and T(($V$)) | 2729 |
| exp | exponential of a tensor series | 2729 |
| log | logarithm of a tensor series | 2729 |
| $\pi_{\text{sym}}$ | the "unormalized" projection to symmetric tensors | 2729 |
| | Posets and lattices | |
| $[n]$ | the set $\{1, \ldots, n\}$ | 2733 |
| $(P, \leq)$ | a poset is a set $P$ with a partial order given by $\leq$ | 2733 |
| $C \subset P$ | a totally ordered subset $C$ of $P$ is called a chain | 2733 |
| $C \subset P$ | a subset $C$ of $P$ with no elements comparable is called an antichain | 2733 |
| $m$ | $m : P \times P \to \mathbb{R}$ the Möbius function associated to the poset $P$ | 2733 |
| | Posets of (ordered) partitions | |
| $\mathcal{P}(S)$ | poset of partitions of a finite set $S$ | 2734, 2735 |
| $\mathcal{P}(d)$ | poset of partitions of $S := [d] := \{1, \ldots, d\}$ | 2734 |
| Orp($P$) | poset of ordered partitions of a finite poset $P$ | 2735 |
| $P_{(n_1, \ldots, n_k)}$ | poset $P$ of disjoint unions of $k$ mutually non-comparable chains $C_i$ of length $n_i$ | 2736 |
| $\mathbf{a}!$ | "factorial" of an ordered partition $\mathbf{a}$ counts the functions that represent the partition $\mathbf{a}$ | 2736 |
| $\mathcal{A}(\mathbf{a})$ | antichain ancestry of an ordered partition $\mathbf{a}$ | 2737 |
| $\mathfrak{d}(\mathbf{a})$ | "height" of an ordered partition $\mathbf{a}$ | 2737 |
| | Signature moments and cumulants | |
| $\int dx^{\otimes m}$ | shorthand for $\int_{0 \leq t_1 \leq \cdots \leq t_m \leq T} dx(t_1) \otimes \cdots \otimes dx(t_m)$ | 2730 |
| $S(x) = (S_m(x))_{m \geq 0}$ | the signature of continuous bounded variation path $x$, $S_m(x) := \int dx^{\otimes m}$ | 2730 |
| $\mathbf{x}$ | geometric $p$-rough path | 2740 |
| $\boldsymbol{\mu}_{\mathbf{X}} = (\mu_{\mathbf{X}}^m)_{m \geq 0}$ | signature moments $\mathbb{E}[\mathbf{X}_{0,T}]$ of a random geometric $p$-rough path $\mathbf{X}$ | 2741 |
| $\boldsymbol{\kappa}_{\mathbf{X}} = (\kappa_{\mathbf{X}}^m)_{m \geq 0}$ | signature cumulants $\log \mathbb{E}[\mathbf{X}_{0,T}]$ of a random geometric $p$-rough path $\mathbf{X}$ | 2741 |
| $\mu_{\mathbf{T}}(\mathbf{a})$ | generalized moments of a T(($V$))-valued random variable $\mathbf{T}$ | 2741 |
| $\kappa_{\mathbf{T}}(\mathbf{a})$ | generalized cumulants of a T(($V$))-valued random variable $\mathbf{T}$ | 2741 |
| | Estimators for signature cumulants | |
| $\hat{k}_n$ | Tukey's polykay; min. variance unbiased estimator for classical cumulants | 2748 |
| $\hat{\boldsymbol{\kappa}}_n$ | signature polykay; min. variance unbiased estimator for signature cumulants | 2749 |

## 2. The lattice of (ordered) partitions

Given a finite set $S$, the set of partitions $\mathcal{P}(S)$ of $S$ is a classical combinatorial object. In this section, we show that if $S$ is replaced by a set $P$ that has an additional partial order structure, then a natural generalisation of $\mathcal{P}(S)$ are the "ordered partitions" $\mathrm{Orp}(P)$. Analogous to the role the poset structure of $\mathcal{P}(S)$ for $S = \{1, \ldots, d\}$ plays for cumulants of $\mathbb{R}^d$-valued random variables, Theorem 1.1, we will see in Section 3 that the poset structure of $\mathrm{Orp}(P)$ for $P = \{\tau_1\} \cup \cdots \cup \{\tau_k\}$ with $\tau_1, \ldots, \tau_k$ tuples formed from $\{1, \ldots, d\}$, plays a similar important role for our signature cumulants, Theorem 1.1 and Theorem 3.6.

We briefly recall partially ordered sets (henceforth posets), lattices and some important examples such as partitions (classical cumulants), and non-crossing partitions (free cumulants) in Section 2.1 and proceed to introduce the lattice of ordered partitions and study its properties in Section 2.2.

### 2.1. Posets, lattices, and Möbius inversion

Throughout this paper, we use the notation $[n]$ for the set $\{1, \ldots, n\}$ where $n \in \mathbb{N}$. Both $[n]$ and $\mathbb{N}$ are always considered with their normal total order "$\leq$". For a function $f : S \to U$ between sets, we denote by $\ker(f) := \{f^{-1}(u) : u \in U\}$.

A *partially ordered set* (poset) is a pair $(P, \leq)$ where $P$ is a set and "$\leq$" is a partial order on $P$. It is customary to refer to just $P$ as the poset. A function $f : P_1 \to P_2$ between two posets is said to be *order-preserving* if $f(x) \geq f(y)$ in $P_2$ whenever $x \geq y$ in $P_1$, the set of such functions is denoted $\mathrm{Hom}(P_1, P_2)$. We say that a subset $C \subseteq P$ is a *chain* if it is totally ordered and an *antichain* if no element of $C$ is comparable to any other element of $C$.

An important function on any poset $P$ with finite intervals is the *Möbius function*, defined recursively as:

$$m(\mathbf{a}, \mathbf{b}) = \begin{cases} 1 & \text{if } \mathbf{a} = \mathbf{b}, \\ -\sum_{\mathbf{a} \leq \mathbf{c} < \mathbf{b}} m(\mathbf{a}, \mathbf{c}) & \text{if } \mathbf{a} < \mathbf{b}, \\ 0 & \text{otherwise.} \end{cases}$$

The most important application is that for finite $P$, and $f : P \to \mathbb{C}$ and $F : P \to \mathbb{C}$ related by

$$F(\mathbf{a}) = \sum_{\mathbf{b} \leq \mathbf{a}} f(\mathbf{b}),$$

the function $f$ can be recovered from $F$ by "Möbius inversion"

$$f(\mathbf{a}) = \sum_{\mathbf{b} \leq \mathbf{a}} m(\mathbf{b}, \mathbf{a}) F(\mathbf{b}).$$

If every two elements of $(P, \leq)$ have a unique minimal upper bound and unique maximal lower bound then we say that $(P, \leq)$ is a *lattice*. We refer the reader to [23] for many more results and examples of posets and lattices. Below we recall two lattices that motivate signature cumulants and ordered partitions.

*The lattice of partitions*

A *partition* **a** of a finite set $S$ is a set of disjoint subsets $a_1, \ldots, a_k$ of $S$ such that their union equals $S$, or equivalently $\mathbf{a} = \ker(f)$ for some $f : S \to \mathbb{N}$. We call the subsets $a_1, \ldots, a_k$ the *blocks* of the partition $\mathbf{a} = \{a_1, \ldots, a_k\}$ and denote by $|\mathbf{a}|$ the number of blocks of **a**. We denote by $\mathcal{P}(S)$ the *set of partitions of* $S$ and make it into a poset by endowing it with the following partial order: for $\mathbf{a}, \mathbf{b} \in \mathcal{P}(S)$

$$\mathbf{a} \leq \mathbf{b}$$

holds if every block of **a** is contained in some block of **b**. In this case, we say that **a** is *finer* than **b**. It is well known that the pair $(\mathcal{P}(S), \leq)$ is a lattice. We will typically denote $\mathcal{P}([d])$ simply by $\mathcal{P}(d)$ and if $\mathbf{a} \in \mathcal{P}(S)$ and $U \subseteq S$ then we use the notation $\mathbf{a} \cap U := \{a \cap U \mid a \in \mathbf{a}\} \in \mathcal{P}(U)$. As pointed out in [19], a direct application of Möbius inversion yields a generalisation of the bijection between cumulants and moments.

**Proposition 2.1 ([19]).** *For a $\mathbb{R}^d$-valued random variable $X = (X_1, \ldots, X_d)$ define the generalised moments and cumulants of a partition $\mathbf{a} = \{a_1, \ldots, a_k\} \in \mathcal{P}(d)$ as*

$$\mu_X(\mathbf{a}) = \prod_{i=1}^k \mathbb{E}(X^{a_i}), \qquad \kappa_X(\mathbf{a}) = \prod_{i=1}^k \kappa(X_{a_i}), \tag{2.1}$$

*where we denote*

$$X^a = \prod_{i \in a} X_i, \quad X_a = \bigtimes_{i \in a} X_i$$

*for a set $a \subseteq [d]$. Then*

$$\mu_X(\mathbf{a}) = \sum_{\mathbf{b} \leq \mathbf{a}} \kappa_X(\mathbf{b}), \qquad \kappa_X(\mathbf{a}) = \sum_{\mathbf{b} \leq \mathbf{a}} \mu_X(\mathbf{b}) m(\mathbf{b}, \mathbf{a}).$$

In the special case where **a** is the partition with only one block, Proposition 2.1 reduces to the well-known formula

$$\kappa(X_1, \ldots, X_d) = \sum_{\mathbf{a} \in \mathcal{P}(d)} (-1)^{|\mathbf{a}|-1} (|\mathbf{a}| - 1)! \mu_X(\mathbf{a})$$

which is Item 1 of Theorem 1.1. We point the reader to [19,20] for many more results that can be derived with this approach.

*The lattice of non-crossing partitions*

If a set $S$ is totally ordered, then a partition $\mathbf{a} \in \mathcal{P}(S)$ is said to be *crossing* if there exists $x_1 < x_2 < x_3 < x_4$ in $S$ such that $x_1, x_3$ and $x_2, x_3$ belong to two distinct blocks of **a**. If one requires that partitions be compatible with the order on $S$ in sense of not being crossing, then one arrives at the set of *non-crossing partitions* of $S$, typically denoted by NC($S$). When endowed with the partial order inherited from $\mathcal{P}(S)$, NC($S$) becomes a lattice. This lattice naturally appears in free probability where, if one requires that Proposition 2.1 holds, one arrives at the *free cumulants*.

## 2.2. The lattice of ordered partitions

We now introduce a lattice with a richer structure by introducing *ordered partitions*. In spite of the fact that they are defined in a natural way, the literature on them is limited. To the best of the authors knowledge, they first appear in the paper [24] by T. Sturm, and a similar notion appears in the PhD-thesis of R. Stanley [22]. More recently, they have found use in combinatorial algebraic topology, [12]. In this section, we derive some basic results that we will need in the sequel; in Section 3, we see that they are naturally linked to signature cumulants.

Recall that set of partitions of a finite set $S$ can be defined as follows:

**Definition 2.1.** Let $S$ be a finite set. We call

$$\mathcal{P}(S) := \{\ker(f) \mid f : S \to \mathbb{N}\}$$

the set of partions of $S$.

In view of this, the following definition is natural.

**Definition 2.2.** Let $(P, \leq)$ be a finite poset. We call

$$\mathrm{Orp}(P) := \{\ker(f) \mid f \in \mathrm{Hom}(P, \mathbb{N})\}$$

the set of ordered partions of $P$.

Equipped with the partial order of refinement both $\mathcal{P}(P)$, and its subset $\mathrm{Orp}(P)$, form a lattice [24], Theorem 21, and the inclusion map $\iota : \mathrm{Orp}(P) \to \mathcal{P}(P)$ is order-preserving. Like in the case of non-crossing partitions, $\mathrm{Orp}(P)$ is not a sublattice of $\mathcal{P}(P)$ in general.

**Example 2.1.** See Figure 1 for some simple examples; in addition:

1. Take $P = \{x_1, x_2, y_1, y_2\}$ with $x_1 \leq x_2$ and $y_1 \leq y_2$ and no other relations. Then the partition $\{x_1, y_1\}, \{x_2, y_2\}$ is ordered since the map $f : P \to \{1, 2\}$, $f(x_1) = f(y_1) = 1$, $f(x_2) = f(y_2) = 2$ is order-preserving. The partition $\{x_1, y_2\}, \{x_2, y_1\}$ is not ordered however since if $f$ is order-preserving such that $f(x_1) = f(y_2)$ and $f(y_1) = f(x_2)$ then we would need $f(x_2) \geq f(x_1) = f(y_2) \geq f(y_1) = f(x_2)$, hence $\ker(f) \neq \{\{x_1, y_2\}, \{x_2, y_1\}\}$.

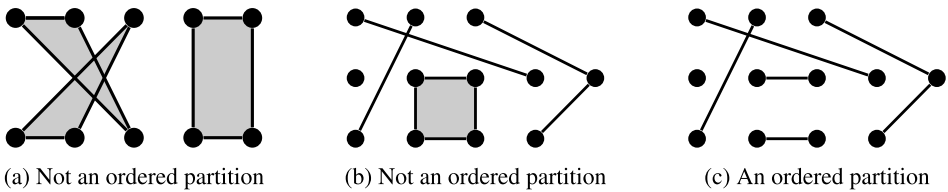

(a) Not an ordered partition  (b) Not an ordered partition  (c) An ordered partition

**Figure 1.** Using the notation from Definition 2.3, one element of $P_{(5,5)}$ and two elements of $P_{(3,5,4)}$ all drawn with the order going left to right. One can check that Figure 1(a), Figure 1(b) are not ordered partitions while Figure 1(c) is.

2. If $P$ is an antichain, then $\mathrm{Orp}(P) = \mathcal{P}(P)$.
3. If $P$ is totally ordered, then $\mathrm{Orp}(P)$ is the set of partitions of $P$ with convex blocks, i.e if $x$ and $y$ belong to the same block then so do all $x \leq z \leq y$.

**Remark 2.1.** Given a partition $\mathbf{a} \in \mathcal{P}(P)$, we say that a sequence $x_0, \ldots x_n$ is an $\mathbf{a}$-*cycle* if $x_0 = x_n$ and for every $i \in \{0, \ldots, n-1\}$ either $x_i < x_{i+1}$ or $x_i, x_{i+1}$ belong to the same block of $\mathbf{a}$. If we say that a partition $\mathbf{a} \in \mathcal{P}(P)$ is $n$-*ordered* if $x_0, \ldots, x_n$ all belong to the same block of $\mathbf{a}$ for any $\mathbf{a}$-cycle $x_0, \ldots, x_n$. Then we get the characterisation that $\mathbf{a}$ is ordered if it is $n$-ordered for every $n \geq 1$ ([25], Lemma 2) and non-crossing if it's 4-ordered.

In this paper, we are mainly concerned with $P$ a disjoint union of mutually non-comparable chains as defined below.

**Definition 2.3.** We denote by $\mathrm{P}_{(n_1, \ldots, n_k)}$ the poset $P = C_1 \cup \cdots \cup C_k$ where each $C_i$ is a chain of length $n_i$ and $x \in C_i$, $y \in C_j$ are comparable if and only if $i = j$.

For a partition $\mathbf{a} \in \mathcal{P}(S)$ the number of blocks is related to the number of functions that represent it via the factorial,

$$|\mathbf{a}|! = \#\{f : S \to [|\mathbf{a}|] \mid \ker(f) = \mathbf{a}\}.$$

This is no longer true for ordered partitions and motivates the following notation.

**Definition 2.4.** For $\mathbf{a} \in \mathrm{Orp}(P)$, we define its *factorial* $\mathbf{a}!$ as

$$\mathbf{a}! := \#\{f \in \mathrm{Hom}(P, [|\mathbf{a}|]) \mid \ker(f) = \mathbf{a}\}.$$

For this next example, recall that if $C \subseteq P$ and $\mathbf{a} \in \mathcal{P}(P)$ we use the notation

$$\mathbf{a} \cap C := \{a \cap C \mid a \in \mathbf{a}\} \in \mathcal{P}(C).$$

**Example 2.2.**

1. If $P = C_1 \cup \cdots \cup C_k = \mathrm{P}_{(|C_1|, \ldots, |C_k|)}$, and $\mathbf{a} = (\mathbf{a} \cap C_1) \cup \cdots \cup (\mathbf{a} \cap C_k)$ is a disjoint union of its intersection with each chain, then $\mathbf{a}!$ is easily seen to be:

$$\mathbf{a}! = \frac{|\mathbf{a}|!}{\prod_{i=1}^{k} |\mathbf{a} \cap C_i|!}.$$

2. If $P$ is an antichain, then $\mathbf{a}! = |\mathbf{a}|!$ for any $\mathbf{a} \in \mathrm{Orp}(P)$.

**Remark 2.2.** Counting the number of order-preserving maps from a poset $P$ to a chain $[n]$ is a classic topic in combinatorics and computer science, where one is typically concerned with bijective maps, so-called *linear extensions*. An important quantity in order theory is the order polynomial $\Omega(n, P)$ and strict order polynomial $\Omega°(n, P)$ of $P$, see [3]. They will not be used here but we remark that for a finite poset $P$ and $n \leq \#P$

$$\sum_{\mathbf{a} \in \mathrm{Orp}(P), |\mathbf{a}| \leq n} \mathbf{a}! = \Omega(n, P), \qquad \sum_{\mathbf{a} \in \mathrm{Orp}(P), |\mathbf{a}| = n} \mathbf{a}! = \Omega°(n, P),$$

which follows directly from the definition.

Another important quantity that appears in the context of signature cumulants is the *antichain ancestry* of $\mathbf{a} \in \mathrm{Orp}(P)$.

**Definition 2.5.** For $\mathbf{a} \in \mathrm{Orp}(P)$, we define the set

$$\mathcal{A}(\mathbf{a}) = \{\mathbf{b} \geq \mathbf{a} : \mathbf{b} \cap C = \mathbf{a} \cap C \text{ for any chain } C \subseteq P\}.$$

We further define $\mathfrak{d}(\mathbf{a})$ as:

$$\mathfrak{d}(\mathbf{a}) = \sum_{\mathbf{b} \in \mathcal{A}(\mathbf{a})} (-1)^{|\mathbf{b}|-1} \frac{\mathbf{b}!}{|\mathbf{b}|}.$$

As it turns out, there is a correspondence between $\mathcal{A}(\mathbf{a})$ and $\{f \in \mathrm{Hom}(P, [|\mathbf{a}|]) \mid \ker(f) = \mathbf{a}\}$, and $\mathfrak{d}(\mathbf{a}) = 0$ whenever $\mathbf{a}$ is minimal in its antichain ancestry. This is spelled out below in the case where $P$ is the disjoint union of exactly two chains.

**Remark 2.3.** If $P$ is an antichain, then this is easy to see since

$$\mathfrak{d}(\mathbf{a}) = \sum_{\mathbf{a} \leq \mathbf{b}} (-1)^{|\mathbf{b}|-1} (|\mathbf{b}| - 1)! = \sum_{\mathbf{a} \leq \mathbf{b} \leq \mathbf{1}} m(\mathbf{b}, \mathbf{1}) = \begin{cases} 1 & \text{if } \mathbf{a} = \mathbf{1}, \\ 0 & \text{otherwise}, \end{cases}$$

where $\mathbf{1}$ is the partition with one block.

**Proposition 2.2.** *If* $P = C_1 \cup C_2 = P_{(|C_1|, |C_2|)}$ *and* $\mathbf{a} \in \mathrm{Orp}(P)$, *then*

$$\mathbf{a}! = \#\mathcal{A}(\mathbf{a}).$$

*Moreover, if* $\mathbf{a} = (\mathbf{a} \cap C_1) \cup (\mathbf{a} \cap C_2)$, *then* $\mathfrak{d}(\mathbf{a}) = 0$.

**Proof.** Assume without loss of generality that $\mathbf{a} \cap C_i = C_i$. Define the set

$$\mathcal{L}(\mathbf{a}) = \{f \in \mathrm{Hom}(P, [|\mathbf{a}|]) \mid \ker(f) = \mathbf{a}\}.$$

If $\{x, y\} \subseteq C_1 \times C_2$ and $\{u, v\} \subseteq C_1 \times C_2$ are two blocks of $\mathbf{a}$, then either $x \leq u$, $y \leq v$ or $u \leq x$, $v \leq y$. Hence, if $f, g \in \mathcal{L}(\mathbf{a})$ and $\{x, y\} \subseteq C_1 \times C_2$ is a block of $\mathbf{a}$, then $f(x) = f(y) = g(x) = g(y)$. So if we let $\mathbf{a}_0^\star, \ldots, \mathbf{a}_k^\star$ be the blocks of $\mathbf{a}$ with exactly two elements and $\mathbf{a}_i$ be the (possibly empty) collection of elements $x_i$ such that $\mathbf{a}_{i-1}^{\star} \cap C_1 < x_i < \mathbf{a}_i^{\star} \cap C_1$ if $x_i \in \mathbf{a}_i \cap C_1$ and $\mathbf{a}_{i-1}^{\star} \cap C_2 < x_i < \mathbf{a}_i^{\star} \cap C_2$ if $x_i \in \mathbf{a}_i \cap C_2$, then we may identify:

$$\mathcal{L}(\mathbf{a}) = \mathcal{L}(\mathbf{a}_1) \times \cdots \times \mathcal{L}(\mathbf{a}_k).$$

Additionally, if $\{x, y\} \subseteq C_1 \times C_2$ is a block of $\mathbf{a}$ and $\mathbf{b} \in \mathcal{A}(\mathbf{a})$ has a block $\{u, v\} \subseteq C_1 \times C_2$, then either:

1. $u = x, v = y$,
2. $u < x, v < y$,
3. $u > x, v > y$.

In view of this, we may make the similar identification:

$$\mathcal{A}(\mathbf{a}) = \mathcal{A}(\mathbf{a}_1) \times \cdots \times \mathcal{A}(\mathbf{a}_k).$$

To show the first claim, we may therefore assume that $\mathbf{a} = C_1 \cup C_2$ so that the blocks of $\mathbf{a}$ are singletons, but this is now clear since in this case, if $\mathbf{b} \in \mathcal{A}(\mathbf{a})$, then $\mathbf{b}$ is the aggregation of some number of blocks of $C_1$ with the same number of blocks in $C_2$. Since $\mathbf{b}$ is ordered, there is exactly one unique way of aggregating them when the blocks are chosen. So there's exactly $\binom{|C_1|}{m}\binom{|C_2|}{m}$ number of ordered partitions $\mathbf{b} \in \mathcal{A}(\mathbf{a})$ that aggregates $m$ blocks of $\mathbf{a}$, hence

$$\#\mathcal{A}(\mathbf{a}) = \sum_{m=0}^{|C_1| \wedge |C_2|} \binom{|C_1|}{m}\binom{|C_2|}{m} = \frac{|C|!}{|C_1|!|C_2|!} = \mathbf{a}!$$

by Remark 1. This shows the first claim. To see the second claim, define the auxiliary sets

$$\mathcal{A}_n(\mathbf{a}) = \{\mathbf{b} \in \mathcal{A}(\mathbf{a}) : |\mathbf{b}| = |\mathbf{a}| - n\}$$

and let $z_n(\mathbf{a}) = \#\mathcal{A}_n(\mathbf{a})$. We know that:

$$\mathbf{a}! = \sum_{n \geq 0} z_n(\mathbf{a}).$$

Noting the identity

$$\sum_{\mathbf{b} \in \mathcal{A}_i(\mathbf{a})} z_j(\mathbf{b}) = \binom{i+j}{i} z_{i+j}(\mathbf{a}),$$

we can write

$$\mathfrak{d}(\mathbf{a}) = \sum_{\mathbf{b} \geq \mathbf{a}} (-1)^{|\mathbf{b}|-1} \frac{\mathbf{b}!}{|\mathbf{b}|}$$

$$= \sum_{n \geq 0} \frac{1}{|\mathbf{a}| - n} (-1)^{|\mathbf{a}|-n-1} \sum_{\mathbf{b} \in \mathcal{A}_n(\mathbf{a})} \sum_{k \geq 0} z_k(\mathbf{b})$$

$$= (-1)^{|\mathbf{a}|-1} \sum_{n\geq 0} \frac{1}{|\mathbf{a}|-n} (-1)^n \sum_{k\geq 0} \binom{n+k}{k} z_{n+k}(\mathbf{a})$$

$$= (-1)^{|\mathbf{a}|-1} \sum_{m\geq 0} z_m(\mathbf{a}) \sum_{n+k=m} \frac{1}{|\mathbf{a}|-n} (-1)^n \binom{m}{k}$$

$$= (-1)^{|\mathbf{a}|-1} \sum_{m\geq 0} (-1)^m \frac{z_m(\mathbf{a})}{(|\mathbf{a}|-m)\binom{|\mathbf{a}|}{m}}.$$

Now fix $\mathbf{a} = (\mathbf{a} \cap C_1) \cup (\mathbf{a} \cap C_2)$ and let $|\mathbf{a} \cap C_1| = m_1$ and $|\mathbf{a} \cap C_2| = m_2$. By the same reasoning as before, we have

$$z_m(\mathbf{a}) = \binom{m_1}{m}\binom{m_2}{m},$$

and we get

$$\mathfrak{d}(\mathbf{a}) = (-1)^{|\mathbf{a}|-1} \sum_{m=0}^{m_1 \wedge m_2} (-1)^m \frac{\binom{m_1}{m}\binom{m_2}{m}}{(m_1 + m_2 - m)\binom{m_1+m_2}{m}}$$

$$= (-1)^{|\mathbf{a}|-1} \frac{1}{m_1 + m_2} \sum_{m=0}^{m_1 \wedge m_2} (-1)^m \frac{\binom{m_1}{m}\binom{m_2}{m}}{\binom{m_1+m_2-1}{m}}$$

$$= (-1)^{|\mathbf{a}|-1} \frac{1}{m_1 + m_2} {}_2F_1(-m_1, -m_2; -m_1 - m_2 + 1, 1) = 0,$$

where ${}_2F_1$ is the Gaussian hypergeometric function. □

## 3. The signature cumulant and its properties

The lattice of ordered partitions from Section 2 provides us with the combinatorial language to prove properties about signature cumulants in analogy with classical cumulants. Section 3.1 recalls geometric rough paths and in Section 3.2 we introduce generalised signature moments and cumulants. Finally, in Section 3.3 we prove Theorem 3.6 that is our "path-space" version of the classical cumulant results, Theorem 1.1.

### 3.1. Geometric rough paths

Given a basis $\{e_1, \ldots, e_d\}$ for a vector space $V$, the set $\{e_\tau := e_{i_1} \otimes \cdots \otimes e_{i_m} : \tau = (i_1, \ldots, i_m) \in [d]^m\}$ is a basis for $V^{\otimes m}$. Throughout, we denote with $[d]^\star = \bigcup_{m\geq 0}[d]^m$ the set of tuples of arbitrary length.

**Definition 3.1.** The shuffle product $\shuffle : T(V) \times T(V) \to T(V)$ is defined by linear extension of the map

$$e_{\tau_1} \shuffle e_{\tau_2} = \sum_{\tau \in Sh(\tau_1, \tau_2)} e_\tau, \quad \tau_1, \tau_2 \in [d]^\star,$$

where the sum is taken over all shuffles[3] of $\tau_1$ and $\tau_2$.

**Remark 3.1.** If $\tau_1 = (i_1, \ldots, i_m)$, $\tau_2 = (i_{m+1}, \ldots, i_{m+n})$, then their shuffles can equivalently be defined as:

$$\text{Sh}(\tau_1, \tau_2) = \left\{ f(\tau_1, \tau_2) = (i_{f^{-1}(1)}, \ldots, i_{f^{-1}(m+n)}) \mid \text{bijections } f \in \text{Hom}\big(P_{(m,n)}, P_{(m+n)}\big) \right\}.$$

**Definition 3.2.** Given $p \geq 1$, a *weakly geometric p-rough path* in $V$ is a map $\mathbf{x} : [0, T]^2 \to T((V))$ such that for all $0 \leq s < t < u \leq T$

1. $\langle \mathbf{x}_{s,t}, f \shuffle g \rangle = \langle \mathbf{x}_{s,t}, f \rangle \langle \mathbf{x}_{s,t}, g \rangle$ for all $f, g \in T(V)$,
2. $\mathbf{x}_{s,t} \otimes \mathbf{x}_{t,u} = \mathbf{x}_{s,u}$,
3. $\|\mathbf{x}\|_{p\text{-var}} < \infty$,

where

$$\|\mathbf{x}\|_{p\text{-var}} = \max_{1 \leq m \leq p} \sup_D \left( \sum_{t_i \in D} |\mathbf{x}_{t_i, t_{i+1}}^m|^{p/m} \right)^{1/p},$$

and the sum is taken over all dissections[4] of $[0, T]$. The space of weakly geometric $p$-rough paths in $V$ is denoted by $\mathcal{G}_p^w(V)$, and $\mathcal{G}_p^w(V)$-valued random variables will be denoted by $\mathbf{X}$.

*Random rough paths*

Our main interest here is that many stochastic processes $X$ in the state space $V$ can be naturally lifted to random (weakly) geometric rough paths $\mathbf{X}$.

**Example 3.1.**

1. If the sample paths of $X$ have finite $p$-variation for some $p < 2$, then the lift may be constructed using Young integration. For instance, any piecewise $C^1$ path has finite 1-variation.
2. If $X$ is a continuous semi-martingale, then the lift may be constructed using Stratonovich integration.
3. If $X$ is a Gaussian process with i.i.d. components and a covariance function of finite $p$-variation for some $p < 4$, then there is a canonical lift of $X$, see [9].

---

[3]A shuffle of tuples $\tau_1 = (i_1, \ldots, i_m) \in [d]^m$ and $\tau_2 = (j_1, \ldots, j_n) \in [d]^n$ is given by mapping $(\sigma_1, \ldots, \sigma_{m+n}) := (i_1, \ldots, i_m, j_1, \ldots, j_n)$ to $(\sigma_{\pi(1)}, \ldots, \sigma_{\pi(m+n)})$ where $\pi$ is a permutation on $[m+n]$ such that $\pi(1) < \cdots \pi(m)$ and $\pi(m+1) < \cdots < \pi(m+n)$.
[4]By a *dissection* of $[0, T]$ we mean a finite collection $t_1 < \cdots < t_{|D|} \in [0, T]$

## 3.2. Signature moments and cumulants

In view of Theorem 1.1, the following definition is not surprising.

**Definition 3.3.** We call
$$\mu_{\mathbf{X}} := (\mu_{\mathbf{X}}^m)_{m \geq 0} := \mathbb{E}[\mathbf{X}_{0,T}] \in T((V))$$
the *signature moments* of $\mathbf{X}$ and we call
$$\kappa_{\mathbf{X}} := (\kappa_{\mathbf{X}}^m)_{m \geq 0} := \log \mathbb{E}[\mathbf{X}_{0,T}] \in T((V))$$
the *signature cumulant* of $\mathbf{X}$ whenever these quantities are well-defined.

It is well known that the signature moments characterise the law of a signature if they decay fast enough, see [5]. Throughout the remainder of Section 3, $\mathbf{X}$ denotes a geometric rough path for which the signature moments are well-defined.

*Generalised signature cumulants*

As in the classical case, a generalisation of moments and cumulants will be useful; see, for example, [19,20] and Proposition 2.1. Applied to tuples that are single indices the definition below recovers the classical generalized moments and cumulants.

**Definition 3.4.** Let $\tau_1, \ldots, \tau_k \in [d]^\star$. We denote with $\mathcal{P}(\tau_1, \ldots, \tau_k)$ the partitions of the set $\{\tau_1\} \cup \cdots \cup \{\tau_k\}$ and $\mathrm{Orp}(\tau_1, \ldots, \tau_k)$ the ordered partitions of the same set endowed with the poset structure of $P_{(|\tau_1|, \ldots, |\tau_k|)}$. For $\mathbf{a} = \{a_1, \ldots, a_n\} \in \mathcal{P}(\tau_1, \ldots, \tau_k)$, we define $a_i^j := \tau_j \cap a_i$.

**Definition 3.5.** Given a random variable $\mathbf{T}$ with values in $T((V))$ and $\mathbf{a} \in \mathcal{P}(\tau_1, \ldots, \tau_k)$, we call
$$\mu_{\mathbf{T}}(\mathbf{a}) = \prod_{i=1}^{|\mathbf{a}|} \mathbb{E}\big[\langle \mathbf{T}, e_{a_i^1}\rangle \cdots \langle \mathbf{T}, e_{a_i^k}\rangle\big]$$
the generalised $\mathbf{a}$-moments of $\mathbf{T}$ and
$$\kappa_{\mathbf{T}}(\mathbf{a}) = \prod_{i=1}^{|\mathbf{a}|} \kappa\big(\langle \mathbf{T}, e_{a_i^1}\rangle, \ldots, \langle \mathbf{T}, e_{a_i^k}\rangle\big)$$
the generalised $\mathbf{a}$-cumulants of $\mathbf{T}$. Here, by convention, $a_i^1, \ldots, a_i^k$ only runs over the non-empty $a_i^j$'s.

**Example 3.2.** If $\mathbf{a} = 1235|67|4 \in \mathcal{P}(1234, 567)$, then
$$\mu_{\mathbf{T}}(\mathbf{a}) = \mathbb{E}\big(\langle \mathbf{T}, e_{123}\rangle \langle \mathbf{T}, e_5\rangle\big) \mathbb{E}\big(\langle \mathbf{T}, e_{67}\rangle\big) \mathbb{E}\big(\langle \mathbf{T}, e_4\rangle\big),$$
$$\kappa_{\mathbf{T}}(\mathbf{a}) = \kappa\big(\langle \mathbf{T}, e_{123}\rangle, \langle \mathbf{T}, e_5\rangle\big) \kappa\big(\langle \mathbf{T}, e_{67}\rangle\big) \kappa\big(\langle \mathbf{T}, e_4\rangle\big).$$

**Lemma 3.1.** *For a random variable* $\mathbf{T}$ *with values in* $T((V))$ *and* $\mathbf{a} \in \mathcal{P}(\tau_1, \ldots, \tau_k)$, *it holds that*

$$\mu_{\mathbf{T}}(\mathbf{a}) = \sum_{\mathbf{b} \leq \mathbf{a}, \mathbf{a} \in \mathcal{A}(\mathbf{b})} \kappa_{\mathbf{T}}(\mathbf{b}).$$

**Proof.** To see this, fix $\mathbf{a} \in \mathcal{P}(\tau_1, \ldots, \tau_k)$ and let $P$ be the poset $a^1 \cup \cdots \cup a^k$ with an antichain ordering. Let $\Phi : \{\tau_1\} \cup \cdots \cup \{\tau_k\} \to P$ be the map such that $\Phi(x) = a^i_j$ if $x \in a^i_j$. Given a random variable $\mathbf{T}$, define the random variable $\mathbf{T}'$ taking values in $\mathbb{R}^P$ by

$$\langle \mathbf{T}', e_i \rangle := \langle \mathbf{T}, e_{\Phi^{-1}(i)} \rangle, \quad i \in P,$$

then we can then write

$$\mu_{\mathbf{T}}(\mathbf{a}) = \mu_{\mathbf{T}'}(\Phi(\mathbf{a})),$$

where $\Phi$ acts on $\mathbf{a}$ by

$$\Phi(\mathbf{a}) = \{\Phi(a) : a \in \mathbf{a}\} \in \mathcal{P}(P).$$

Note that $\mu_{\mathbf{T}'}(\Phi(\mathbf{a}))$ is the moment as defined in Equation (2.1) since $P$ is an antichain. By Proposition 2.1, we get

$$\mu_{\mathbf{T}'}(\Phi(\mathbf{a})) = \sum_{\mathbf{b}' \leq \Phi(\mathbf{a})} \kappa_{\mathbf{T}'}(\mathbf{b}')$$

which, by definition of $\mathbf{T}'$ yields

$$\mu_{\mathbf{T}}(\mathbf{a}) = \sum_{\mathbf{b} = \Phi^{-1}(\mathbf{b}'), \mathbf{b}' \leq \Phi(\mathbf{a})} \kappa_{\mathbf{T}}(\mathbf{b}).$$

So it is enough to note that $\mathbf{b} = \Phi^{-1}(\mathbf{b}')$ with $\mathbf{b}' \leq \Phi(\mathbf{a})$ if and only if $\mathbf{b} \leq \mathbf{a}$ and $\mathbf{a} \in \mathcal{A}(\mathbf{b})$. □

### 3.3. Cumulants, moments, and independence

We show Theorem 1.2, restated for weakly geometric rough paths at the end of this section.

**Lemma 3.2.** *For any $p \geq 1$ and $\mathcal{G}^w_p(V)$-valued random variable $\mathbf{X}$ and any $\tau \in [d]^\star$*

$$\langle \kappa_{\mathbf{X}}, e_\tau \rangle = \sum_{\mathbf{a} \in \mathrm{Orp}(\tau)} \frac{(-1)^{|\mathbf{a}|-1}}{|\mathbf{a}|} \mu_{\mathbf{X}}(\mathbf{a}),$$

$$\langle \mu_{\mathbf{X}}, e_\tau \rangle = \sum_{\mathbf{a} \in \mathrm{Orp}(\tau)} \frac{1}{|\mathbf{a}|!} \kappa_{\mathbf{X}}(\mathbf{a}).$$

**Proof.** Let $\overline{\mu_X} = \mu_X - 1$, that is,

$$\langle \overline{\mu_X}, e_\tau \rangle = \begin{cases} \langle \mu_X, e_\tau \rangle & \text{if } \tau \neq \varnothing, \\ 0 & \text{otherwise.} \end{cases}$$

We may write

$$\log \mu_X = \log(1 + \overline{\mu_X}) = \sum_{n \geq 1} \frac{(-1)^{n-1}}{n} (\overline{\mu_X})^{\otimes n} = \sum_{n \geq 1} \frac{(-1)^{n-1}}{n} \left( \sum_{m \geq 1} \mu_X^m \right)^{\otimes n},$$

$$\mu_X = \exp \kappa_X = 1 + \sum_{n \geq 1} \frac{1}{n!} \kappa_X^{\otimes n} = 1 + \sum_{n \geq 1} \frac{1}{n!} \left( \sum_{m \geq 1} \kappa_X^m \right)^{\otimes n}.$$

By identifying coordinates we see that:

$$\langle \kappa_X, e_\tau \rangle = \sum_{n \geq 1} \sum_{\tau_1 \cdots \tau_n = \tau} \frac{(-1)^{n-1}}{n} \prod_i^n \langle \mu_X, e_{\tau_i} \rangle,$$

$$\langle \mu_X, e_\tau \rangle = \sum_{n \geq 1} \sum_{\tau_1 \cdots \tau_n = \tau} \frac{1}{n!} \prod_i^n \langle \kappa_X, e_{\tau_i} \rangle,$$

where the sum is over all de-concatenations of $\tau$ into non-empty sub-words $\tau_1, \ldots, \tau_k$. The assertion now follows from the identification made in Item 3 of example 2.1. $\square$

**Proposition 3.3.** *For any $p \geq 1$ and $\mathcal{G}_p^w(V)$-valued random variable $\mathbf{X}$ and any tuple $\boldsymbol{\tau} = (\tau_1, \ldots, \tau_k), \tau_1, \ldots, \tau_k \in [d]^\star$[5]*

$$\langle \kappa_X, e_{\tau_1} \shuffle \cdots \shuffle e_{\tau_k} \rangle = \sum_{\mathbf{a} \in \mathrm{Orp}(\boldsymbol{\tau})} (-1)^{|\mathbf{a}|-1} \frac{\mathbf{a}!}{|\mathbf{a}|} \mu_X(\mathbf{a}).$$

**Proof.** Let $N = |\tau_1| + \cdots + |\tau_k|$ and write $|\boldsymbol{\tau}| = (|\tau_1|, \ldots, |\tau_k|)$ for brevity. Define the set

$$\mathrm{Sh}(\boldsymbol{\tau}) = \{ f(\boldsymbol{\tau}) \mid \text{bijective } f \in \mathrm{Hom}(P_{|\boldsymbol{\tau}|}, [N]) \}.$$

By Lemma 3.2 and Remark 3.1 we may write:

$$\langle \kappa_X, e_{\tau_1} \shuffle \cdots \shuffle e_{\tau_k} \rangle = \sum_{w \in \mathrm{Sh}(\boldsymbol{\tau})} \sum_{\mathbf{a} \in \mathrm{Orp}(w)} \frac{(-1)^{|\mathbf{a}|-1}}{|\mathbf{a}|} \mu_X(\mathbf{a}) = \sum_{\mathbf{a} \in \mathrm{Orp}(N)} \frac{(-1)^{|\mathbf{a}|-1}}{|\mathbf{a}|} \sum_{w \in \mathrm{Sh}(\boldsymbol{\tau})} \mu_X(w|\mathbf{a}),$$

---

[5]To avoid having to deal with partitions on multisets we will assume throughout that all indices are distinct. This is treated similarly in [19,20], that is, all indices are treated as distinct even if the same index appears multiple times.

where we write $w|\mathbf{a}$ for the partition of $\tau$ with blocks $(w|\mathbf{a})_i = \bigcup_{j \in a_i} w_j$. So it is enough to show that for any $q \in \mathbb{N}$

$$\sum_{\substack{\mathbf{a} \in \text{Orp}(N) \\ |\mathbf{a}|=q}} \sum_{w \in \text{Sh}(\tau)} \mu_\mathbf{X}(w|\mathbf{a}) = \sum_{\substack{\mathbf{b} \in \text{Orp}(\tau) \\ |\mathbf{b}|=q}} \mathbf{b}! \mu_\mathbf{X}(\mathbf{b}).$$

Fix some $\mathbf{a} = \ker(f_\mathbf{a}) = \{a_1, \ldots, a_q\} \in \text{Orp}(N)$ and some $w \in \text{Sh}(\tau)$. Since $[N]$ is totally ordered $f_\mathbf{a} : [N] \to [q]$ is uniquely determined by $\mathbf{a}$. Since $w \in \text{Sh}(\tau)$ there exists some linear extension $f_w : P_{|\tau|} \to [N]$ so that $w = f_w(\tau)$. Define a map

$$\Phi : \text{Orp}(N) \times \text{Sh}(\tau) \to \text{Hom}(P_{|\tau|}, [q]), \qquad (\mathbf{a}, w) \mapsto f_\mathbf{a} \circ f_w.$$

In addition, for any $\mathbf{a} \in \text{Orp}(N)$, $f_\mathbf{b} \in \text{Hom}(P_{|\tau|}, [q])$, define the set:

$$\mathcal{D}(\mathbf{a}, f_\mathbf{b}) := \{w \in \text{Sh}(\tau) \mid f_\mathbf{b} = f_\mathbf{a} \circ f_w\}$$

and let $H(f_\mathbf{b})$ be the unique element of $\text{Orp}(N)$ such that $\# f_\mathbf{a}^{-1}(i) = \# f_\mathbf{b}^{-1}(i)$ for every $i \in [q]$. Note the following:

1. $\Phi$ is surjective since if $f_\mathbf{b} \in \text{Hom}(P_{|\tau|}, [q])$, then by fixing some linear extensions on every block of $\mathbf{b}$ we may factor $f_\mathbf{b} = \tilde{f}_\mathbf{b} \circ h$ where $h : P_{|\tau|} \to [N]$ is a linear extension of $P_{|\tau|}$. Then $\Phi(\ker(\tilde{f}_\mathbf{b}), h(\tau)) = f_\mathbf{b}$.
2. The kernel of $\Phi$ decomposes as follows:

$$\Phi^{-1}(f_\mathbf{b}) = \{(\mathbf{a}, w) \in \text{Orp}(N) \times \text{Sh}(\tau) \mid \mathbf{a} = H(f_\mathbf{b}), w \in \mathcal{D}(\mathbf{a}, f_\mathbf{b})\}.$$

3. For any $f_\mathbf{b} \in \text{Hom}(P_{|\tau|}, [q])$ such that $\ker(f_\mathbf{b}) = \mathbf{b}$ and $\mathbf{a} = H(f_\mathbf{b})$, item 1 of Definition 3.2 implies that

$$\sum_{w \in \mathcal{D}(\mathbf{a}, f_\mathbf{b})} \mu_\mathbf{X}(w|\mathbf{a}) = \mu_\mathbf{X}(\mathbf{b}).$$

In view of this, we can write

$$\sum_{\substack{\mathbf{a} \in \text{Orp}(N) \\ |\mathbf{a}|=q}} \sum_{w \in \text{Sh}(\tau)} \mu_\mathbf{X}(w|\mathbf{a})$$

$$= \sum_{f_\mathbf{b} \in \text{Hom}(P_{|\tau|}, [q])} \sum_{(\mathbf{a}, w) \in \Phi^{-1}(f_\mathbf{b})} \mu_\mathbf{X}(w|\mathbf{a})$$

$$= \sum_{f_\mathbf{b} \in \text{Hom}(P_{|\tau|}, [q])} \sum_{w \in \mathcal{D}(H(f_\mathbf{b}), \mathbf{b})} \mu_\mathbf{X}(w|H(f_\mathbf{b}))$$

$$= \sum_{f_\mathbf{b} \in \mathrm{Hom}(P_{|\tau|}, [q])} \mu_\mathbf{X}(\mathbf{b})$$

$$= \sum_{\substack{\mathbf{b} \in \mathrm{Orp}(\tau) \\ |\mathbf{b}| = q}} \mathbf{b}! \mu_\mathbf{X}(\mathbf{b}).$$

□

**Corollary 3.4.** *Under the same assumptions as Proposition* 3.3:

$$\langle \kappa_\mathbf{X}, e_{\tau_1} \sqcup \cdots \sqcup e_{\tau_k} \rangle = \sum_{\mathbf{a} \in \mathrm{Orp}(\tau)} \partial(\mathbf{a}) \kappa_\mathbf{X}(\mathbf{a}).$$

**Proof.** Using Remark 3.1, we may write

$$\langle \kappa_X, e_{\tau_1} \sqcup \cdots \sqcup e_{\tau_k} \rangle = \sum_{\mathbf{a} \in \mathrm{Orp}(\tau)} (-1)^{|\mathbf{a}|-1} \frac{\mathbf{a}!}{|\mathbf{a}|} \mu_\mathbf{X}(\mathbf{a})$$

$$= \sum_{\mathbf{a} \in \mathrm{Orp}(\tau)} (-1)^{|\mathbf{a}|-1} \frac{\mathbf{a}!}{|\mathbf{a}|} \sum_{\mathbf{b} \leq \mathbf{a}, \mathbf{a} \in \mathcal{A}(\mathbf{b})} \kappa_\mathbf{X}(\mathbf{b})$$

$$= \sum_{\mathbf{b} \in \mathrm{Orp}(\tau)} \kappa_\mathbf{X}(\mathbf{b}) \sum_{\mathbf{a} \in \mathcal{A}(\mathbf{b})} (-1)^{|\mathbf{a}|-1} \frac{\mathbf{a}!}{|\mathbf{a}|}$$

$$= \sum_{\mathbf{b} \in \mathrm{Orp}(\tau)} \partial(\mathbf{b}) \kappa_\mathbf{X}(\mathbf{b}).$$

□

**Example 3.3.** If $\tau = (i_1, \ldots, i_m) \in [d]^\star$, then by Remark 2.3 we can conclude that:

$$\langle \pi_{\mathrm{sym}}(\kappa_\mathbf{X}), e_\tau \rangle = \langle \kappa_\mathbf{X}, e_{i_1} \sqcup \cdots \sqcup e_{i_m} \rangle = \kappa(\langle \mathbf{X}_{0,T}, e_{i_1} \rangle, \ldots, \langle \mathbf{X}_{0,T}, e_{i_m} \rangle),$$

and one retrieves the classical cumulant on the increment $\mathbf{X}_{0,T}$.

**Proposition 3.5.** *Let $I, J \subseteq [d]$ be two sets of coordinates and let $\mathbf{X}$ take values in $\mathcal{G}_p^w(V)$ for some $p \geq 1$. Then*

1. $\mathbb{E}[\langle \mathbf{X}_{0,T}, e_{\tau_1} \rangle \langle \mathbf{X}_{0,T}, e_{\tau_2} \rangle] = \langle \mu_\mathbf{X}, e_{\tau_1} \rangle \langle \mu_\mathbf{X}, e_{\tau_2} \rangle, \forall \tau_1 \in I^\star, \tau_2 \in J^\star$
   *if and only if*
2. $\langle \kappa_\mathbf{X}, e_{\tau_1} \sqcup e_{\tau_2} \rangle = 0, \forall \tau_1 \in I^\star, \tau_2 \in J^\star$ *non-empty*,

*where $I^\star := \{(i_1, \ldots, i_m) \mid i_1, \ldots, i_m \in I, m \geq 1\}$.*

**Proof.** For any fixed $\tau_1 \in I^\star, \tau_2 \in J^\star$ we say that $\mathbf{a} \in \mathrm{Orp}(\tau_1, \tau_2)$ is non-degenerate if $\mathbf{a} \neq (\mathbf{a} \cap \tau_1) \cup (\mathbf{a} \cap \tau_2)$. By Corollary 3.4 and Proposition 2.2 we may write

$$\langle \kappa_\mathbf{X}, e_{\tau_1} \sqcup e_{\tau_2} \rangle = \kappa(\langle \mathbf{X}_{0,T}, e_{\tau_1} \rangle, \langle \mathbf{X}_{0,T}, e_{\tau_2} \rangle) + \sum_{\substack{\mathbf{a} \in \mathrm{Orp}(\tau_1, \tau_2), |\mathbf{a}| \geq 2 \\ \mathbf{a} \text{ non-degenerate}}} \partial(\mathbf{a}) \kappa_\mathbf{X}(\mathbf{a}).$$

Noting that 1 is equivalent to $\kappa(\langle \mathbf{X}_{0,T}, e_{\tau_1}\rangle, \langle \mathbf{X}_{0,T}, e_{\tau_2}\rangle) = 0$ for every $\tau_1 \in I^\star$, $\tau_2 \in J^\star$ non-empty, the assertion follows by induction on the length of $\tau_1, \tau_2$. □

**Remark 3.2.** We would like to thank the referee for pointing out the following algebraic proof of Proposition 3.5 which relies on the bialgebra structure of $T((V))$ instead of the combinatorial results established here.

Write $\mathbb{R}^I$, $\mathbb{R}^J$ for the vector spaces spanned by $\{e_i \mid i \in I\}$ and $\{e_j \mid j \in J\}$ respectively and denote by $\pi_I : \mathbb{R}^d \to \mathbb{R}^I$, $\pi_J : \mathbb{R}^d \to \mathbb{R}^J$ the associated projection maps. These naturally extend to projection maps $\pi_I^\star : T((\mathbb{R}^d)) \to T((\mathbb{R}^I))$, $\pi_J^\star : T((\mathbb{R}^d)) \to T((\mathbb{R}^J))$ which are both graded algebra morphisms. Denoting by $\Delta$ the map defined by algebraic extension of $\Delta v = 1 \otimes v + v \otimes 1$, it is well known that $\Delta$ is a graded algebra morphism and is dual to ⧢ in the sense that

$$\langle x, f \shuffle g\rangle = \langle \Delta x, f \otimes g\rangle$$

for $x \in T((V))$, $f, g \in T(V^\star)$, see, for example, [18], Proposition 1.9. Hence, Item 2 is equivalent to

$$\langle \Delta \kappa_\mathbf{X}, e_{\tau_1} \otimes e_{\tau_2}\rangle = \langle \kappa_\mathbf{X} \otimes 1 + 1 \otimes \kappa_\mathbf{X}, e_{\tau_1} \otimes e_{\tau_2}\rangle$$

for any $\tau_1 \in I^\star$, $\tau_2 \in J^\star$; The above can be restated as

$$(\pi_I^\star \otimes \pi_J^\star)\Delta \kappa_\mathbf{X} = (\pi_I^\star \otimes \pi_J^\star)(\kappa_\mathbf{X} \otimes 1 + 1 \otimes \kappa_\mathbf{X}). \quad (3.1)$$

Since $\pi_I^\star \otimes \pi_J^\star : T((\mathbb{R}^d)) \otimes T((\mathbb{R}^d)) \to T((\mathbb{R}^I)) \otimes T((\mathbb{R}^J))$ and $\Delta : T((\mathbb{R}^d)) \to T((\mathbb{R}^d)) \otimes T((\mathbb{R}^d))$ are graded algebra morphisms they factor through the exponential map, which we may apply to both sides of Equation (3.1) to get

$$(\pi_I^\star \otimes \pi_J^\star)\Delta e^{\kappa_\mathbf{X}} = (\pi_I^\star \otimes \pi_J^\star)e^{\kappa_\mathbf{X} \otimes 1 + 1 \otimes \kappa_\mathbf{X}} = (\pi_I^\star \otimes \pi_J^\star)(e^{\kappa_\mathbf{X}} \otimes e^{\kappa_\mathbf{X}})$$

since $e^{\kappa_\mathbf{X}} = \mu_\mathbf{X}$ by definition, this is equivalent to

$$(\pi_I^\star \otimes \pi_J^\star)\Delta \mu_\mathbf{X} = (\pi_I^\star \otimes \pi_J^\star)(\mu_\mathbf{X} \otimes \mu_\mathbf{X})$$

which is equivalent to Item 1.

To summarise this section, we have now proven Theorem 1.2. To state it appropriately, we first need to define the space of *unparametrised* weakly geometric rough paths.

**Definition 3.6.** The space of *unparametrised* weakly geometric p-rough paths is defined as the quotient space

$$\mathcal{R}_p^w(V) = \mathcal{G}_p^w(V)/\sim,$$

where $\sim$ denotes so called *tree-like equivalence* of rough paths and is essentially factoring out different ways of parametrising it. See Appendix B for details.

**Remark 3.3.** It is well known in the literature that $\mathbf{X} \mapsto \mathbf{X}_{0,T}$ is a bijection for $\mathbf{X} \in \mathcal{R}_p^w(V)$ [1], Theorem 1.1.

Signature cumulants, ordered partitions, and independence of stochastic processes 2747

Putting everything together and using the notation $\mathbf{X}|_I = (\langle \mathbf{X}, e_\tau \rangle)_{\tau \in I^\star}$ for a set of coordinates $I \subseteq [d]$ shows the following result.

**Theorem 3.6.** *Let $p \geq 1$ and $\mathbf{X}$ be a $\mathcal{R}_p^w(V)$-valued random variable such that $\mathbb{E}[|\langle \mathbf{X}_{0,T}, f \rangle|] < \infty$ for every $f \in T(V)$. Then for $\tau, \tau_1, \ldots, \tau_k \in [d]^\star$ it holds that*

1. (compensated moments)

$$\langle \kappa_{\mathbf{X}}, e_{\tau_1} \shuffle \cdots \shuffle e_{\tau_k} \rangle = \sum_{\mathbf{a} \in \mathrm{Orp}(\tau_1, \ldots, \tau_k)} (-1)^{|\mathbf{a}|-1} \frac{\mathbf{a}!}{|\mathbf{a}|} \mu_{\mathbf{X}}(\mathbf{a}),$$

2. (moment bijection)

$$\langle \mu_{\mathbf{X}}, e_\tau \rangle = \sum_{\mathbf{a} \in \mathrm{Orp}(\tau)} \frac{1}{|\mathbf{a}|!} \kappa_{\mathbf{X}}(\mathbf{a}),$$

3. (independence) *If $I, J \subseteq [d]$ and the joint law $(\mathbf{X}|_I, \mathbf{X}|_J)$ is characterised by*

$$\big\{ \mathbb{E}\big[\langle \mathbf{X}_{0,T}, e_{\tau_1} \rangle \langle \mathbf{X}_{0,T}, e_{\tau_2} \rangle\big] \mid \tau_1 \in I^\star, \tau_2 \in J^\star \big\},$$

*then $\mathbf{X}|_I$ and $\mathbf{X}|_J$ are independent if and only if*

$$\langle \kappa_{\mathbf{X}}, e_{\tau_1} \shuffle e_{\tau_2} \rangle = 0 \quad \text{for any } \tau_1 \in I^\star, \tau_2 \in J^\star.$$

A sufficient condition for Item 3 is that the individual signature moments decay sufficiently fast, see Appendix A. This can be verified for many models; an alternative is to use "normalised signatures" [6], for which the normalised signature moments always characterise the law which is useful for "black box" approaches as used in machine learning.

## 4. Estimating signature cumulants

To apply the previous results in practice, one has to take into account that only a finite number of sample paths are available and needs to derive efficient estimators for signature cumulants. For classical cumulants, Tukey introduce so-called *Polykays* in [26]. We recall these in Section 4.1 and present them in our formulation. In Section 4.2, we then introduce such polykays for signature cumulants and show their optimality and other properties as estimators in Proposition 4.2. Finally, we apply these results to give an application to drift and volatility estimation for SDEs in Example 4.3.

### 4.1. Tukey's polykays

Given a partition $\mathbf{a} = \{a_1, \ldots, a_k\}$ of $[d]$, the *symmetric means* are polynomials in $n \times d$ variables $(x_i^1, \ldots, x_i^d)_{i=1}^n$ defined as follows:

$$\hat{\mu}_n(\mathbf{a}) = \frac{1}{(n)_k} \sum_{1 \leq i_1, i_2, \ldots, i_k \leq n}^{\neq} x_{i_1}^{a_1} \cdots x_{i_k}^{a_k},$$

where $(n)_k = n(n-1)\cdots(n-k)$, $x_i^S = \prod_{j \in S} x_i^j$ and the notation $\sum^{\neq}$ means that the sum is taken over all non-equal indices. The *Polykays* – originally introduced by Tukey in [26] – are defined as

$$\hat{k}_n(\mathbf{a}) = \sum_{\mathbf{b} \leq \mathbf{a}} m(\mathbf{b}, \mathbf{a}) \hat{\mu}_n(\mathbf{b}),$$

where $m$ is the Möbius function of the partition lattice. See, for example, [14,20] for more on polykays.

**Example 4.1.**

$$\hat{\mu}_n(12) = \frac{1}{n} \sum_{i=1}^n x_i^1 x_i^2, \qquad \hat{\mu}_n(1|2) = \frac{1}{n(n-1)} \sum_{1 \leq j \neq i \leq n} x_i^1 x_j^2.$$

For $\mathbf{a} \in \mathcal{P}(m)$ and i.i.d. random variables $(X_i)_{i=1}^n = (X_i^1, \ldots, X_i^m)_{i=1}^n$ the symmetric means and polykays are unbiased estimators of moments and cumulants respectively, that is,

$$\mathbb{E}\hat{\mu}_n(\mathbf{a}) = \mu_X(\mathbf{a}), \qquad \mathbb{E}\hat{k}_n(\mathbf{a}) = \kappa_X(\mathbf{a}).$$

**Lemma 4.1.** *For any finite disjoint sets $S_1$, $S_2$ and $\mathbf{a}_1 \in \mathcal{P}(S_1)$, $\mathbf{a}_2 \in \mathcal{P}(S_2)$:*

$$\mathrm{Cov}(\hat{k}_n(\mathbf{a}_1), \hat{k}_n(\mathbf{a}_2)) = \frac{1}{n} \sum_{\mathbf{c} \in V(\mathbf{a}_1, \mathbf{a}_2)} \sum_{\substack{\mathbf{b} \leq \mathbf{c} \\ \mathbf{b} \nmid \mathbf{a}_1, \mathbf{b} \nmid \mathbf{a}_2}} \kappa_X(\mathbf{b}) + \mathcal{O}(n^{-2}),$$

*where*

$$V(\mathbf{a}_1, \mathbf{a}_2) = \{\mathbf{c} \in \mathcal{P}(S_1 \cup S_2) : \mathbf{c} \cap S_1 = \mathbf{a}_1, \mathbf{c} \cap S_2 = \mathbf{a}_2, |\mathbf{c}| = |\mathbf{a}_1| + |\mathbf{a}_2| - 1\},$$

*and the notation $\mathbf{b} \nmid \mathbf{a}$ means that no block of $\mathbf{b}$ is a proper subset of a block of $\mathbf{a}$.*

**Proof.** It follows from [20], Proposition 3.1, that for any $\mathbf{a}_1 \in \mathcal{P}(S_1)$, $\mathbf{a}_2 \in \mathcal{P}(S_2)$

$$\mathbb{E}(\hat{k}_n(\mathbf{a}_1)\hat{k}_n(\mathbf{a}_2)) = \sum_{\mathbf{b} \in \mathcal{P}(S_1 \cup S_2)} c(\mathbf{b}) \kappa_X(\mathbf{b}),$$

$$c(\mathbf{b}) = \sum_{\mathbf{c} \geq \mathbf{b}} \frac{(n)_\mathbf{c}}{(n)_{\mathbf{c}_1}(n)_{\mathbf{c}_2}} m(\mathbf{c}_1, \mathbf{a}_1) m(\mathbf{c}_2, \mathbf{a}_2),$$

where $(n)_\mathbf{c} = n(n-1)\cdots(n-|\mathbf{c}|)$ and $\mathbf{c}_i = \mathbf{c} \cap S_i$. By [20], Proposition 3.2, we know that $c(\mathbf{b}) = 0$ if any block of $\mathbf{b}$ is a proper subset of a block of $\mathbf{a}_1$ or $\mathbf{a}_2$, hence the only $\mathcal{O}(1)$ term is $\kappa_X(\mathbf{a}_1)\kappa_X(\mathbf{a}_2)$.

Unless there exists some $\mathbf{c} \geq \mathbf{b}$ such that $|\mathbf{c}| = |\mathbf{c}_1| + |\mathbf{c}_2| - 1$, then $c(\mathbf{b}) \in \mathcal{O}(n^{-2})$ and $\frac{(n)_\mathbf{c}}{(n)_{\mathbf{c}_1}(n)_{\mathbf{c}_2}} = \frac{1}{n} + \mathcal{O}(n^{-2})$ if this is the case. Additionally, since no block of $\mathbf{b}$ is a proper subset of a block of $\mathbf{a}_1$, $\mathbf{a}_2$ this must also be true for such a $\mathbf{c} \geq \mathbf{b}$. Combined with the constraint $|\mathbf{c}| = |\mathbf{c}_1| + |\mathbf{c}_2| - 1$ this implies that $\mathbf{c}_1 = \mathbf{a}_1$, $\mathbf{c}_2 = \mathbf{a}_2$, so $\mathbf{c} \in V(\mathbf{a}_1, \mathbf{a}_2)$. In view of this, we may write:

$$\mathrm{Cov}(\hat{k}_n(\mathbf{a}_1), \hat{k}_n(\mathbf{a}_2)) = \frac{1}{n} \sum_{\substack{\mathbf{b} \in \mathcal{P}(S_1 \cup S_2) \\ \mathbf{b} \nmid \mathbf{a}_1, \mathbf{b} \nmid \mathbf{a}_2}} \kappa_X(\mathbf{b}) \sum_{\mathbf{c} \geq \mathbf{b}, \mathbf{c} \in V(\mathbf{a}_1, \mathbf{a}_2)} m(\mathbf{c}_1, \mathbf{a}_1) m(\mathbf{c}_2, \mathbf{a}_2) + \mathcal{O}(n^{-2})$$

$$= \frac{1}{n} \sum_{\mathbf{c} \in V(\mathbf{a}_1, \mathbf{a}_2)} \sum_{\substack{\mathbf{b} \leq \mathbf{c} \\ \mathbf{b} \nmid \mathbf{a}_1, \mathbf{b} \nmid \mathbf{a}_2}} \kappa_X(\mathbf{b}) + \mathcal{O}(n^{-2}).$$

$\square$

**Remark 4.1.** If one is interested in computing the magnitude of the error in the $\mathcal{O}(n^{-2})$ term then one would use the formula from [20],

$$\mathrm{Cov}(\hat{k}_n(\mathbf{a}_1)\hat{k}_n(\mathbf{a}_2)) = \sum_{\mathbf{b}} c(\mathbf{b}) \kappa_X(\mathbf{b}), \qquad c(\mathbf{b}) = \sum_{\mathbf{c} \geq \mathbf{b}} \frac{(n)_{\mathbf{c}}}{(n)_{\mathbf{c}_1}(n)_{\mathbf{c}_2}} m(\mathbf{c}_1, \mathbf{a}_1) m(\mathbf{c}_2, \mathbf{a}_2)$$

and sum over all $\mathbf{b} \in \mathcal{P}(S_1 \cup S_2)$ such that $\frac{(n)_{\mathbf{c}}}{(n)_{\mathbf{c}_1}(n)_{\mathbf{c}_2}} \in \mathcal{O}(n^{-2})$ for all $\mathbf{c} \geq \mathbf{b}$.

**Example 4.2.** If $\mathbf{a} \in \mathcal{P}(S)$ is a partition with blocks of size at most 2, then Lemma 4.1 shows that:

$$\mathrm{Var}(\hat{k}_n(\mathbf{a})) = \frac{1}{n} \sum_{a_i, a_j \in \mathbf{a}} \left( \kappa_X(a_i a_j) + \kappa_X(a_i^1 a_j^1) \kappa_X(a_i^2 a_j^2) + \kappa_X(a_i^1 a_j^2) \kappa_X(a_i^2 a_j^1) \right) \prod_{k \neq i,j} \kappa_X(a_k)^2$$
$$+ \mathcal{O}(n^{-2}).$$

**Definition 4.1.** If $\mathbf{X}_1, \ldots, \mathbf{X}_n \sim \mathbf{X}$ is an i.i.d. sequence of random variables in $\mathcal{G}_p^w(V)$ and $\boldsymbol{\tau} = (\tau_1, \ldots, \tau_n)$ then we define

$$\hat{\boldsymbol{\kappa}}_n(\boldsymbol{\tau}) = \sum_{\mathbf{a} \in \mathrm{Orp}(\boldsymbol{\tau})} (-1)^{|\mathbf{a}|-1} \frac{\mathbf{a}!}{|\mathbf{a}|} \hat{\mu}_n(\mathbf{a}) = \sum_{\mathbf{a} \in \mathrm{Orp}(\boldsymbol{\tau})} \mathfrak{d}(\mathbf{a}) \hat{k}_n(\mathbf{a}).$$

### 4.2. Signature polykays

The following proposition describes properties of $\hat{\boldsymbol{\kappa}}_n(\boldsymbol{\tau})$. Recall that for a function $f$ of $k$ variables and i.i.d. samples $X_1, \ldots, X_n \sim X$, the *U-statistic* of $f(X)$ is defined as

$$U(f)_n = \mathbb{E}_{\sigma \in \mathrm{Sym}(n)} [f(X_{\sigma(1)}, \ldots, X_{\sigma(k)})],$$

where the expectation is taken over the uniform measure of permutations $\sigma : [n] \to [n]$, see [27].

**Proposition 4.2.** *Let $\mathbf{X}_1, \ldots, \mathbf{X}_n \sim \mathbf{X}$ be an i.i.d. sequence such that $\mathbb{E}|\langle \mathbf{X}, e_\tau \rangle| < \infty$ for every $\tau \in [d]^\star$, then for any $\boldsymbol{\tau} = (\tau_1, \ldots, \tau_k)$:*

1. $\mathbb{E}\hat{\boldsymbol{\kappa}}_n(\boldsymbol{\tau}) = \langle \kappa_{\mathbf{X}}, e_{\tau_1} \shuffle \cdots \shuffle e_{\tau_k}\rangle$,
2. $\hat{\boldsymbol{\kappa}}_n(\boldsymbol{\tau})$ is minimum variance in the family of unbiased polynomial estimators of $\langle \kappa_{\mathbf{X}}, e_{\tau_1} \shuffle \cdots \shuffle e_{\tau_k}\rangle$,
3. $\hat{\boldsymbol{\kappa}}_n(\boldsymbol{\tau}) \to \langle \kappa_{\mathbf{X}}, e_{\tau_1} \shuffle \cdots \shuffle e_{\tau_k}\rangle$ a.s. and in $L^p$ for every $p \in [1, \infty)$ as $n \to \infty$,
4. *For any finite collection $S$ of index tuples, the family $(\hat{\boldsymbol{\kappa}}_n(\boldsymbol{\tau}))_{\boldsymbol{\tau} \in S}$ is asymptotically normal with asymptotic covariance* $\mathcal{V}(\boldsymbol{\tau}_1, \boldsymbol{\tau}_2) = \sum_{\substack{\mathbf{a}_1 \in \mathrm{Orp}(\boldsymbol{\tau}_1) \\ \mathbf{a}_2 \in \mathrm{Orp}(\boldsymbol{\tau}_2)}} \mathfrak{d}(\mathbf{a}_1)\mathfrak{d}(\mathbf{a}_2)\mathcal{V}(\mathbf{a}_1, \mathbf{a}_2)$, where $\mathcal{V}(\mathbf{a}_1, \mathbf{a}_2) = \sum_{\mathbf{c} \in V(\mathbf{a}_1, \mathbf{a}_2)} \sum_{\substack{\mathbf{b} \leq \mathbf{c} \\ \mathbf{b} \nmid \mathbf{a}_1, \mathbf{b} \nmid \mathbf{a}_2}} \kappa_{\mathbf{X}}(\mathbf{b})$.

**Proof.** U-statistics are minimum variance in the family of unbiased polynomial estimators, converging a.s. and in $L^p$, as well as asymptotically normal [27], Section 12. Moreover, if $U, V$ are U-statistics for $f(X), g(X)$, then $U + V$ and $(U, V)$ are U-statistics for $(f+g)(X)$ and $(f, g)(X)$ respectively. Since $\hat{k}_n(\mathbf{a})$ is a U-statistic for $\kappa_X(\mathbf{a})$, $\hat{\boldsymbol{\kappa}}_n(\boldsymbol{\tau})$ is a U-statistic for $\langle \kappa_{\mathbf{X}}, e_{\tau_1} \shuffle \cdots \shuffle e_{\tau_k}\rangle$ and this extends to any finite collection. It only remains to show the asymptotic variance. By Lemma 4.1, we get

$$\mathbb{E}\hat{\boldsymbol{\kappa}}_n(\boldsymbol{\tau}_1)\hat{\boldsymbol{\kappa}}_n(\boldsymbol{\tau}_2) - \mathbb{E}\hat{\boldsymbol{\kappa}}_n(\boldsymbol{\tau}_1)\mathbb{E}\hat{\boldsymbol{\kappa}}_n(\boldsymbol{\tau}_2)$$
$$= \sum_{\substack{\mathbf{a}_1 \in \mathrm{Orp}(\boldsymbol{\tau}_1) \\ \mathbf{a}_2 \in \mathrm{Orp}(\boldsymbol{\tau}_2)}} \mathfrak{d}(\mathbf{a}_1)\mathfrak{d}(\mathbf{a}_2)\big(\mathbb{E}\big[\hat{k}_n(\mathbf{a}_1)\hat{k}_n(\mathbf{a}_2)\big] - \kappa(\mathbf{a}_1)\kappa(\mathbf{a}_2)\big)$$
$$= \sum_{\substack{\mathbf{a}_1 \in \mathrm{Orp}(\boldsymbol{\tau}_1) \\ \mathbf{a}_2 \in \mathrm{Orp}(\boldsymbol{\tau}_2)}} \mathfrak{d}(\mathbf{a}_1)\mathfrak{d}(\mathbf{a}_2)\mathrm{Cov}\big(\hat{k}_n(\mathbf{a}_1), \hat{k}_n(\mathbf{a}_2)\big)$$
$$= \frac{1}{n}\sum_{\substack{\mathbf{a}_1 \in \mathrm{Orp}(\boldsymbol{\tau}_1) \\ \mathbf{a}_2 \in \mathrm{Orp}(\boldsymbol{\tau}_2)}} \mathfrak{d}(\mathbf{a}_1)\mathfrak{d}(\mathbf{a}_2)\mathcal{V}(\mathbf{a}_1, \mathbf{a}_2) + \mathcal{O}(n^{-2}).$$
□

**Example 4.3 (Estimating a diffusion with constant drift and volatility).** Denote by $B$ a standard $d$-dimensional Browian motion and consider the process

$$X_t = bt + \sigma B_t,$$

where $b \in \mathbb{R}^d$ and $\sigma \in (\mathbb{R}^d)^{\otimes 2}$. It is well known that $X$ can be lifted to a geometric rough path $\mathbf{X}$ and that

$$\mu_{\mathbf{X}} := \mathbb{E}[\mathbf{X}_{0,1}] = \exp\left(b + \frac{1}{2}\sigma\right)$$

see for instance [9], Exercise 3.22. Hence, the first two signature cumulants and signature moments of $\mathbf{X}$ are

$$\mu_{\mathbf{X}}^1 = b, \qquad\qquad \kappa_{\mathbf{X}}^1 = b,$$
$$\mu_{\mathbf{X}}^2 = \frac{1}{2}(b^{\otimes 2} + \sigma), \qquad \kappa_{\mathbf{X}}^2 = \frac{1}{2}\sigma.$$

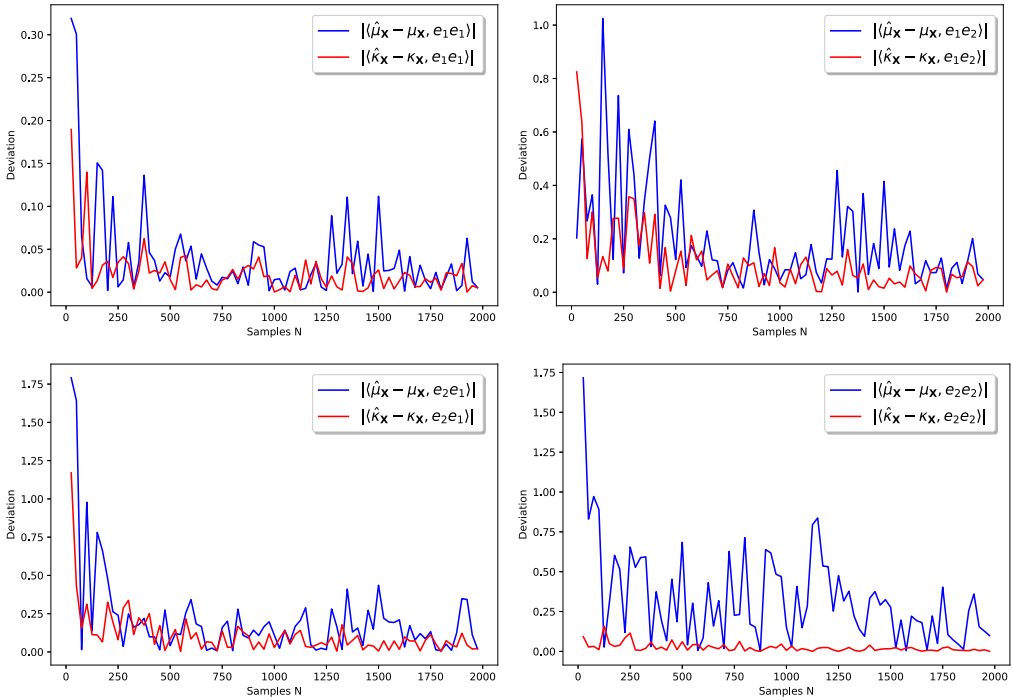

**Figure 2.** We simulated a 2-dimensional diffusion $X$ with drift $b = \binom{1}{10}$ and volatility $\sigma = I_2$ the identity matrix. The plots show the absolute difference $|\hat{\mu}_{\mathbf{X}}^2 - \mu_{\mathbf{X}}^2|$, and $|\hat{\kappa}_{\mathbf{X}}^2 - \kappa_{\mathbf{X}}^2|$, between the theoretical values of the second signature moment (blue) and the second signature cumulant (red) and their estimators as the number of samples $N$ varies in the range $25 \leq N \leq 2000$.

By Item 2 of Theorem 3.6, both $\mu_{\mathbf{X}}$ and $\kappa_{\mathbf{X}}$ characterise the law of the process $X$, but if one wants to learn the law of the process $X$ from observing sample trajectories of $X$, then it is much more efficient to work with signature cumulants than with signature moments, thus echoing a classic insight from statistics for vector-valued data, see Example 1.1. To see this, assume we are given an observation of the process $(X_t(\omega))_{t \in [0,N]}$ over a time interval $[0, N]$ for a large $N \in \mathbb{N}$. We divide this interval into $N$ pieces of length 1 and calculate the signature of $X$ over each of these shorter intervals. Since $X$ is strong Markov, this results in $N$ independent samples $\mathbf{X}_1, \ldots, \mathbf{X}_N$ of the signature $\mathbf{X}$ of $(X_t)_{t \in [0,1]}$. The unbiased estimators introduced above are given as

$$\hat{\mu}_{\mathbf{X}}^1 = \frac{1}{N} \sum_{i=1}^{N} \mathbf{X}_i^1, \qquad \hat{\kappa}_{\mathbf{X}}^1 = \frac{1}{N} \sum_{i=1}^{N} \mathbf{X}_i^1,$$

$$\hat{\mu}_{\mathbf{X}}^2 = \frac{1}{N} \sum_{i=1}^{N} \mathbf{X}_i^2, \qquad \hat{\kappa}_{\mathbf{X}}^2 = \frac{1}{N} \sum_{i=1}^{N} \mathbf{X}_i^2 - \frac{1}{2N(N-1)} \sum_{1 \leq i \neq j \leq N} \mathbf{X}_i^1 \otimes \mathbf{X}_j^1.$$

We now show that the second signature moment typically has a higher variance than the second signature cumulant. In fact, a direct calculation detailed in Appendix C shows that

$$\text{Var}(\langle \hat{\boldsymbol{\mu}}_\mathbf{X}, e_i e_j \rangle) = \text{Var}(\langle \hat{\boldsymbol{\kappa}}_\mathbf{X}, e_i e_j \rangle) + c_{ij},$$

where the last term $c_{ij}$ is explicitly computed. For example, if $\sigma = I_d$ it reduces to

$$c_{ij} = \frac{1}{N} \begin{cases} \frac{1}{2}(b_i)^4 + 2(b_i)^2 - \frac{1}{2(N-1)} & \text{if } i = j, \\ \frac{1}{2}(b_i)^2(b_j)^2 + \frac{1}{4}(b_i)^2 + \frac{1}{4}(b_j)^2 - \frac{1}{4(N-1)} & \text{otherwise.} \end{cases}$$

and we see that the second signature cumulant estimator has a lower variance than the second signature moment estimator whenever either $b$ or $N$ is sufficiently large, see Figure 2 for a simulated comparison.

## Appendix A: Independence of rough paths on $\mathbb{R}^d$

For an element $x \in T((V))$, let $\pi_n : T((V)) \to V^{\otimes n}$ be the projection map and denote by $x^n = \pi_n x$. Further define:

$$E(V) = \left\{ x \in T((V)) : \sum_{n \geq 0} \|x^n\| \lambda^n < \infty, \forall \lambda > 0 \right\},$$

$$G(V) = \left\{ g \in E : \Delta(g) = g \otimes g, g \neq 0 \right\},$$

where $\Delta$ is the map defined by extension of $\Delta(v) = v \otimes 1 + 1 \otimes v$ for $v \in V$. It follows from [13], Section 2.2.5, that signatures of weakly geometric rough paths in $V$ take values in $G(V)$. Endow $E(V)$ with the locally convex topology induced by the family of norms $\|x\|_\lambda = \sum_{n \geq 0} \|x^n\| \lambda^n$, $\lambda > 0$ and $G(V)$ with the subspace topology. It follows from [5], Sections 2 and 3, that:

1. $G(V)$ is Polish,
2. If $f \in E^\star$ there exists some $\lambda > 0$ such that $\|f \circ \pi_n\| \leq \lambda^n$ and $\sum_{n=0}^N f \circ \pi_n \to f$ pointwise,
3. If $X$ is a random variable that takes values in $G(V)$ and $\mathbb{E}X \in E(V)$, then for $f \in E^\star$, $\mathbb{E}f(X) = \sum_{n \geq 0} \langle \mathbb{E}X^n, f \rangle$.

By [5], Section 4, there exists a subset of $E(V)^\star$ such that when restricted to $G(V)$ forms a $\star$-subalgebra of $C_b(G(V))$ that separates the points. Denote this family of functions by $\mathcal{C}(G)$.

We briefly recall the Stone-Weierstrass theorem for Radon measures ([2], Exercise 7.14.79): If $\mathcal{X}$ is a topological space and $\mu, \nu$ are two Radon measures on $\mathcal{X}$. Then if $\mathcal{F} \subseteq C_b(\mathcal{X})$ is an algebra of functions that separates the points of $\mathcal{X}$, then $\mu = \nu$ is and only if $\mu(f) = \nu(f)$ for every $f \in \mathcal{F}$.

The next lemma is a slight generalisation of [5], Proposition 6.1, and its proof closely resembles that of [5], Proposition 3.2.

**Lemma A.1.** *If $X_1, \ldots, X_k$ are random variables taking values in $G(V)$ that satisfy*:

$$\sum_{n \geq 0} \|\mathbb{E} X_i^n\| \lambda^n < \infty, \quad \forall \lambda > 0$$

*for $i = 1, \ldots, k$. Then the joint law of $(X_1, \ldots, X_k)$ is uniquely determined by*:

$$\mathbb{E}[\langle X_1, e_{\tau_1} \rangle \cdots \langle X_k, e_{\tau_k} \rangle], \quad \tau_1, \ldots, \tau_k \in [d]^\star. \tag{A.1}$$

**Proof.** The set of finite linear combinations of maps of the form

$$(x_1, \ldots, x_k) \mapsto f_1(x_1) \cdots f_k(x_k), \quad f_1, \ldots, f_k \in \mathcal{C}(G)$$

is a $\star$-subalgebra of $C_b(G^k)$ that separates the points of $G^k$. Since $V$ is Polish, so is $G(V)$ and $G(V)^k$, so by [2], Theorem 7.1.7, the law of $(X_1, \ldots, X_k)$ is Radon and hence by Stone-Weierestrass for Radon measures the assertion will then follow if we can show that Equation (A.1) determines $\mathbb{E}[f_1(X_1) \cdots f_k(X_k)]$ for every $f_1, \ldots, f_k \in \mathcal{C}(G)$.

Fix some $f_1, \ldots, f_k \in \mathcal{C}(G)$ and let $f_i^n = f_i \circ \pi_n$. For every $n \geq 0$, $f_i^n \in T(V)$ so $\mathbb{E} f_i^n(X_i)$ is determined by Equation (A.1). Since we know that $\sum_{n=0}^N f_i^n \to f_i$ pointwise and $\mathbb{E} f_i(X_i) = \sum_{n \geq 0} \mathbb{E} f_i^n(X_i)$ it is enough to show that $\sum_{n_1, \ldots, n_k \geq 0} \mathbb{E}(|f_1^{n_1}(X_1) \cdots f_k^{n_k}(X_k)|) < \infty$ since one may then apply dominated convergence to get:

$$\mathbb{E}(f_1(X_1) \cdots f_k(X_k)) = \mathbb{E}\left(\lim_{N \to \infty} \left(\sum_{n_1=0}^N f_1^{n_1}(X_1)\right) \cdots \left(\sum_{n_k=0}^N f_k^{n_k}(X_k)\right)\right)$$

$$= \sum_{n_1, \ldots, n_k \geq 0} \mathbb{E}(f_1^{n_1}(X_1) \cdots f_k^{n_k}(X_k)).$$

But for any $f \in E^*$ and $X \in G$ it holds that $f(X)^2 = f^{\otimes 2}(X \otimes X) = f^{\otimes 2}(\Delta X)$. Hence $f^2 = f^{\otimes 2} \circ \Delta \in E^*$ on $G$, so

$$\sum_{n_1, \ldots, n_k \geq 0} \mathbb{E}(|f_1^{n_1}(X_1) \cdots f_k^{n_k}(X_k)|) \leq \sum_{n_1, \ldots, n_k \geq 0} \mathbb{E}|f_1^{n_1}(X_1)|^2 + \cdots + \mathbb{E}|f_k^{n_k}(X_k)|^2$$

$$= \sum_{n \geq 0} (f_1^n)^{\otimes 2} \Delta(\mathbb{E} X_1^{2n}) + \cdots + (f_k^n)^{\otimes 2} \Delta(\mathbb{E} X_k^{2n}).$$

By continuity of $f_1, \ldots, f_k, \Delta$ we may pick some $\lambda > 0$ such that $\|(f_i^n)^{\otimes 2} \circ \Delta\| \leq \lambda^{2n}$ for every $i = 1, \ldots, k$. Now the assertion follows since

$$\sum_{n_1, \ldots, n_k \geq 0} \mathbb{E}(|f_1^{n_1}(X_1) \cdots f_k^{n_k}(X_k)|) \leq \sum_{n \geq 0} \lambda^{2n} [\|\mathbb{E} X_1^{2n}\| + \cdots + \|\mathbb{E} X_k^{2n}\|] < \infty.$$

$\square$

For the next proposition, denote by $\mathcal{G}_p(V)$ the closure of $C^1$-paths in $\mathcal{G}_p^w(V)$. Even though $\mathcal{G}_p(V)$ is strictly smaller than $\mathcal{G}_p^w(V)$, one always has the inclusion $\mathcal{G}_p^w(V) \subseteq \mathcal{G}_q(V)$ if $q > p$ [9], Exercise 2.15. Denote by $\mathcal{R}_p(V)$ the image of $\mathcal{G}_p(V)$ in $\mathcal{R}_p^w(V)$.

**Proposition A.2.** *For any $p \geq 1$ and $\mathbf{X}$, $\mathbf{Y}$ taking values in $\mathcal{R}_p(V)$ such that*

$$\sum_{n \geq 0} \|\mathbb{E}\mathbf{X}_{0,T}^n\|\lambda^n < \infty, \quad \sum_{n \geq 0} \|\mathbb{E}\mathbf{Y}_{0,T}^n\|\lambda^n < \infty, \quad \forall \lambda > 0,$$

*then $\mathbf{X}$ and $\mathbf{Y}$ are independent if and only if*

$$\mathbb{E}\big[\langle \mathbf{X}_{0,T}, e_{\tau_1}\rangle \langle \mathbf{Y}_{0,T}, e_{\tau_2}\rangle\big] = \mathbb{E}\big[\langle \mathbf{X}_{0,T}, e_{\tau_1}\rangle\big]\mathbb{E}\big[\langle \mathbf{Y}_{0,T}, e_{\tau_2}\rangle\big], \quad \forall \tau_1, \tau_2 \in [d]^\star.$$

**Proof.** $\mathcal{G}_p(V)$ is separable [9], Exercise 2.12, hence so is $\mathcal{R}_p(V)$ as a continuous image of $\mathcal{G}_p(V)$. So the laws of $\mathbf{X}, \mathbf{Y}$ are radon [2], Theorem 7.1.7, and the map $\mathbf{X} \mapsto \mathbf{X}_{0,T}$ is a continuous injection [1], Theorem 1.1. Since continuous injections into Hausdorff spaces are injective on Radon measures it is enough to show that $\mathbf{X}_{0,T}, \mathbf{Y}_{0,T}$ are independent.

But $\mathbf{X}_{0,T}$ and $\mathbf{Y}_{0,T}$ take values in $G(V)$, so the assertion follows by applying Lemma A.1 to the measures on $G(V) \times G(V)$ defined by $\mu_1(\mathcal{A} \times \mathcal{B}) = \mathbb{P}(\mathbf{X}_{0,T} \in \mathcal{A}, \mathbf{Y}_{0,T} \in \mathcal{B})$, $\mu_2(\mathcal{A} \times \mathcal{B}) = \mathbb{P}(\mathbf{X}_{0,T} \in \mathcal{A})\mathbb{P}(\mathbf{Y}_{0,T} \in \mathcal{B})$ for Borel set $\mathcal{A}, \mathcal{B}$. □

## Appendix B: Tree-like equivalence of paths

For a path $x : [0, T] \to E$, where $E$ is some topological space, denote by $\overleftarrow{x}$ its *time-reversal*, defined as

$$\overleftarrow{x}(t) := x(T - t).$$

For two paths $x : [0, T] \to E$, $y : [0, S] \to E$, denote by $x \star y : [0, T + S] \to E$ their *concatenation*, defined as

$$(x \star y)_t := \begin{cases} x(t) & \text{for } 0 \leq t < T, \\ y(s) + x(T) & \text{for } T \leq t \leq T + S. \end{cases}$$

**Definition B.1 ([1]).** A continuous path $x : [0, T] \to E$ is said to be *tree-like* if there exists an $\mathbb{R}$-tree $\tau$, a continuous map $\varphi : [0, T] \to \tau$ and a map $\psi : \tau \to E$ such that $\varphi(0) = \varphi(T)$ and $x = \psi \circ \varphi$.

With the above definition in mind, we say that two paths $x, y$ are *tree-like equivalent* if $x \star \overleftarrow{y}$ is tree-like. If $\mathbf{x}, \mathbf{y} \in \mathcal{G}_p^w(V)$ then we say that they are tree-like equivalent if the two paths $t \mapsto \mathbf{x}_{0,t}$ and $s \mapsto \mathbf{y}_{0,s}$ are. This induces an equivalence relation $\sim$ on $\mathcal{G}_p^w(V)$ and the corresponding quotient space

$$\mathcal{R}_p^w(V) = \mathcal{G}_p^w(V)/\sim$$

is called the space of *unparametrised* weakly geometric p-rough paths. This space is equipped with the topology induced by the map $\mathcal{R}_p^w(V) \to E(V)$, $\mathbf{x} \mapsto \mathbf{x}_{0,T}$.

**Remark B.1.** Tree-like equivalence is essentially factoring out different ways of parametrising paths. If one wants to make a statement about a stochastic process $X$ and not its unparametrised counterpart, then one would simply look at the rough path lift of $(x_t, t)_{t \geq 0}$ instead since for any two paths $x$, $y$ with the same starting value $(x_t, t) \sim (y_t, t)$ if and only if $x = y$.

**Remark B.2.** Since a rough path $\mathbf{x}$ is a function of increments $(s, t) \mapsto \mathbf{x}_{s,t}$, the lift of a process $x$ does not in general depend on its starting value $x_0$. If one wishes to make a statement about the whole process, then one would lift $(x_t, tx_0)$.

## Appendix C: Details for Example 4.3

For $i, j \in [d]$ denote with $\hat{\mu}_\mathbf{X}(i, j)$ the moment estimator associated to the partition $\mathbf{a} = i|j$ as seen in Example 4.1, then

$$\mathrm{Var}(\langle \hat{\boldsymbol{\mu}}_\mathbf{X}, e_i e_j \rangle) = \mathrm{Var}(\langle \hat{\boldsymbol{\kappa}}_\mathbf{X}, e_i e_j \rangle) - \frac{1}{4}\mathrm{Var}(\hat{\mu}_\mathbf{X}(i,j)) + \mathrm{Cov}(\langle \hat{\boldsymbol{\mu}}_\mathbf{X}, e_i e_j \rangle, \hat{\mu}_\mathbf{X}(i,j)).$$

By using the formula in Remark 4.1, we can compute the formal expressions for the covariances of the polykays associated to the partitions $\mathbf{a} = a_1|a_2$ and $\mathbf{b} = a_3$,

$$\mathrm{Cov}(\hat{k}(a_1|a_2), \hat{k}(a_3)) = \frac{1}{N}\bigl[\kappa(a_1|a_2 a_3) + \kappa(a_2|a_1 a_3) + 2\kappa(a_1|a_2|a_3)\bigr],$$

$$\mathrm{Cov}(\hat{k}(a_1|a_2), \hat{k}(a_1|a_2)) = \frac{1}{N}\bigl[\kappa(a_1|a_1|a_2 a_2) + \kappa(a_2|a_2|a_1 a_1) + 2\kappa(a_1|a_2|a_1 a_2)\bigr]$$

$$+ \frac{1}{N(N-1)}\bigl[\kappa(a_1 a_1|a_2 a_2) + \kappa(a_1 a_2|a_1 a_2)\bigr].$$

With this one sees that

$$\mathrm{Cov}(\langle \hat{\boldsymbol{\mu}}_\mathbf{X}, e_i e_j \rangle, \hat{\mu}_\mathbf{X}(i,j))$$
$$= \frac{1}{2N}\bigl(b_i b_j(\sigma_{ij} + \sigma_{ji}) + b_i b_i \sigma_{jj} + b_j b_j \sigma_{ii} + 2b_i b_j \sigma_{ij} + b_i b_i b_j b_j\bigr),$$

$$\mathrm{Var}(\hat{\mu}_\mathbf{X}(i,j))$$
$$= \frac{1}{N}\bigl(b_i b_j(\sigma_{ij} + \sigma_{ji}) + b_i b_i \sigma_{jj} + b_j b_j \sigma_{ii}\bigr) + \frac{1}{N(N-1)}(\sigma_{ii}\sigma_{jj} + \sigma_{ij}\sigma_{ij}),$$

and hence

$$c_{ij} := \frac{1}{4N}\bigl(2b_i b_i b_j b_j + b_i b_j(\sigma_{ij} + \sigma_{ji}) + b_i b_i \sigma_{jj} + b_j b_j \sigma_{ii} + 4b_i b_j \sigma_{ij}\bigr)$$
$$- \frac{1}{4N(N-1)}(\sigma_{ii}\sigma_{jj} + \sigma_{ij}\sigma_{ij}).$$

## Acknowledgements

PB is supported by the Engineering and Physical Sciences Research Council [EP/R513295/1]. HO is supported by the Engineering and Physical Sciences Research Council [EP/S026347/1] (DATASIG), the Turing Institute, and the Oxford-Man Institute of Quantitative Finance.